\documentclass[graybox,envcountsame,envcountsect]{svmult}

\usepackage{mathptmx}       % selects Times Roman as basic font
\usepackage{helvet}         % selects Helvetica as sans-serif font
\usepackage{courier}        % selects Courier as typewriter font
\usepackage{type1cm}        % activate if the above 3 fonts are
\usepackage{makeidx}         % allows index generation
\usepackage{graphicx}        % standard LaTeX graphics tool
\usepackage{multicol}        % used for the two-column index
\usepackage[bottom]{footmisc}% places footnotes at page bottom
\usepackage{url}

\usepackage{amsfonts}
\usepackage{amsmath}
\usepackage{amssymb}
\usepackage[arrow,matrix]{xy}

\smartqed

\makeindex             % used for the subject index
\DeclareMathOperator{\modu}{mod}
\DeclareMathOperator{\ind}{ind}

\DeclareMathOperator{\Hom}{Hom}
\DeclareMathOperator{\End}{End}

\DeclareMathOperator{\add}{add}
\DeclareMathOperator{\ann}{ann}

\DeclareMathOperator{\rep}{rep}
\DeclareMathOperator{\op}{op}
\DeclareMathOperator{\ext}{ext}
\def\soc{\text{\rm soc}}
\def\rad{\text{\rm rad}}
\def\Tr{\text{\rm Tr}\,}
\def\supp{\text{\rm supp}\,}

\def\gldim{\text{\rm gl\,dim}\,}
\def\pd{\text{\rm pd}}
\def\id{\text{\rm id}}
\def\Ext{\text{\rm Ext}}
\def\Tor{\text{\rm Tor}}
\def\vdim{\text{\rm\bf dim}\,}
\def\gl{\text{\rm gl}\,}
\def\GL{\text{\rm GL}\,}

\newcommand{\BA}{{\mathbb A}}
\newcommand{\BD}{{\mathbb D}}

\begin{document}

\title*{Cycle-finite module categories}
\titlerunning{Cycle-finite module categories}
\author{Piotr Malicki, Jos\'e Antonio de la Pe\~na and Andrzej Skowro\'nski}
\authorrunning{P. Malicki et al.}
\institute{Piotr Malicki \at Faculty of Mathematics and Computer Science, Nicolaus Copernicus University, Chopina 12/18, 87-100 Toru\'n, Poland,
\email{pmalicki@mat.uni.torun.pl}
\and Jos\'e A. de la Pe\~na \at Centro de Investigaci\'on en Mathem\'aticas (CIMAT), Guanajuato, M\'exico,
\email{jap@cimat.mx}
\and Andrzej Skowro\'nski \at Faculty of Mathematics and Computer Science, Nicolaus Copernicus University, Chopina 12/18, 87-100 Toru\'n, Poland,
\email{skowron@mat.uni.torun.pl}}

\maketitle

\abstract*{We describe the structure of module categories of finite dimensional algebras over an algebraically closed field for which
the cycles of nonzero nonisomorphisms between indecomposable finite dimensional modules are finite (do not belong to the infinite Jacobson
radical of the module category). Moreover, geometric and homological properties of these
module categories are exhibited.}

\abstract{We describe the structure of module categories of finite dimensional algebras over an algebraically closed field for which
the cycles of nonzero nonisomorphisms between indecomposable finite dimensional modules are finite (do not belong to the infinite Jacobson
radical of the module category). Moreover, geometric and homological properties of these
module categories are exhibited.}

\section{Introduction}\label{sec:1}

Throughout the article $K$ denotes a fixed algebraically closed field.

By an algebra we mean an associative finite dimensional $K$-algebra with an identity which we shall assume (without loss of generality) to be basic and
connected. For an algebra $A$, by an $A$-module we mean a finite dimensional right $A$-module. We shall denote by $\modu A$ the category of $A$-modules,
by $\ind A$ its full subcategory formed by the indecomposable modules, by $\Gamma_A$ the Auslander-Reiten quiver of $A$, and by $\tau_A$ the
Auslander-Reiten translation $D\Tr$ in $\Gamma_A$. We shall identify an indecomposable $A$-module with the vertex of $\Gamma_A$ corresponding to it.
From Drozd's Tame and Wild Theorem \cite{[Dro]} (see also \cite{[CB1]}) the class of algebras may be divided into two classes. One class consists of
the wild algebras whose representation theory comprises the representation theories of all algebras over $K$ (see \cite[Chapter XIX]{[SS2]}).
The second class consists of the tame algebras for which the indecomposable modules occur, in each dimension $d$, in a finite number of discrete and
a finite number of one-parameter families. Hence, a classification of the finite dimensional modules is only feasible for tame algebras.
It has been shown by Crawley-Boevey \cite{[CB1]} that, if $A$ is a tame algebra, then, for any dimension $d\geq 1$, all but finitely many isomorphism
classes of indecomposable $A$-modules of dimension $d$ are invariant on the action of $\tau_A$, and hence, by a result due to Hoshino \cite{[Ho]}, lie
in stable tubes of rank one in $\Gamma_A$. The indecomposable modules over tame algebras which do not lie in stable tubes of rank one are called discrete.
A distinguished class of tame algebras is formed by the algebras of finite representation type, having only finitely many isomorphism classes of
indecomposable modules, for which the representation theory is presently rather well understood (see \cite{[BGRS]}, \cite{[Bo4]}, \cite{[Bo5]},
\cite{[BoG]}, \cite{[BrG]}). On the other hand, the representation theory of arbitrary tame algebras is still only emerging.
At present the most accessible seem to be
the (tame) algebras of polynomial growth, for which there exists an integer $m$ such that the number of one-parameter families of indecomposable modules
is bounded, in each dimension $d$, by $d^m$. This class of algebras has been subject of intensive research over the last 30 years.

A prominent role in the representation theory of algebras is played by cycles of modules, or more generally cycles of complexes of modules.
Recall that a cycle in a module category $\modu A$ is a sequence
\[ X_0 \buildrel {f_1}\over {\hbox to 6mm{\rightarrowfill}} X_1 \to \cdots \to X_{r-1} \buildrel {f_r}\over {\hbox to 6mm{\rightarrowfill}} X_r=X_0 \]
of nonzero nonisomorphisms in $\ind A$, and the cycle is said to be finite if the homomorphisms $f_1,\ldots, f_r$ do not belong to the infinite Jacobson
radical of $\modu A$.
Following Ringel \cite{[Ri]} a module in $\ind A$ which does not lie on cycle in $\ind A$ is called directing. It has been proved independently
by Peng and Xio \cite{[PX]} and the third named author \cite{[Sk6]} that the Auslander-Reiten quiver $\Gamma_A$ of an arbitrary algebra $A$ contains at most
finitely many $\tau_A$-orbits containing directing modules. Hence, in order to obtain information on nondirecting indecomposable modules of a module
category, we may study properties of cycles in $\modu A$ containing these modules. We also note that, by a result of Ringel \cite{[Ri]} the
support algebras of directing modules are tilted algebras.
Following \cite{[AS4]} an algebra $A$ is said to be cycle-finite if all cycles in $\modu A$ are finite. It has been proved by the
third named author in \cite{[Sk8]} that every cycle-finite algebra $A$ is of polynomial growth and the support algebras of the one-parametric families
of indecomposable $A$-modules are tame concealed algebras (preprojective tilts of the path algebras of Euclidean quivers) and Ringel's tubular algebras,
which are distinguished classes of cycle-finite algebras. The class of cycle-finite algebras is wide and contains the algebras of finite representation
type, the tame tilted algebras \cite{[Ke]}, the tame double tilted algebras \cite{[RS1]}, the tame generalized double tilted algebras \cite{[RS2]},
the tubular algebras \cite{[Ri]}, the iterated tubular algebras \cite{[PTo]}, the tame quasi-tilted algebras \cite{[LS]}, \cite{[Sk10]},
the tame coil and multicoil algebras \cite{[AS4]}, \cite{[AS5]}, \cite{[AS6]}, the tame generalized multicoil algebras \cite{[MS2]}, and the strongly simply connected
algebras of polynomial growth \cite{[Sk9]}. It has been also proved in \cite{[AS1]}, \cite{[AS2]}, \cite{[AS3]} that the class of algebras $A$ for which
the derived category $D^b(\modu A)$ of bounded complexes of $A$-modules is cycle-finite coincides with the class of piecewise hereditary algebras of
Dynkin, Euclidean, and tubular type, and consequently these algebras are also cycle-finite. Moreover, frequently an algebra $A$ admits a Galois covering
$R\to R/G = A$ where $R$ is a cycle-finite locally bounded category and $G$ is an admissible group of automorphisms of $R$, which allows to reduce the
representation theory of $A$ to the representation theory of cycle-finite algebras being finite convex subcategories of $R$. For example, every
selfinjective algebra $A$ of polynomial growth admits a canonical standard form $\overline A$ (geometric socle deformation of $A$) such that $\overline A$
has a Galois covering $R\to R/G = \overline A$, where $R$ is a cycle-finite selfinjective locally bounded category and $G$ is an admissible infinite
cyclic group of automorphisms of $R$, the Auslander-Reiten quiver $\Gamma_{\overline A}$ of $\overline A$ is the orbit quiver $\Gamma_R/G$ of $\Gamma_R$,
and the stable Auslander-Reiten quivers of $A$ and $\overline A$ are isomorphic (see \cite{[Sk1]}, \cite{[Sk12]} for details).
We also mention that, by the main result of \cite{[PS3]}, every algebra $A$ which admits a cycle-finite Galois covering $R\to R/G = A$ with $G$ torsion-free
is tame.

One of the objectives of this article is to describe the structure of the category $\ind A$ of an arbitrary cycle-finite algebra $A$, by showing that it
can be covered by the categories of indecomposable modules of tame generalized multicoil algebras and tame generalized double tilted algebras.
Here, a crucial role will be played by description of support algebras of cyclic components of the Auslander-Reiten quivers of cycle-finite algebras.
The second objective of the article is to exhibit geometric and homological properties of indecomposable modules over cycle-finite algebras.
We are interested in the class of coherent cycle-finite algebras for which all cyclic components of the Auslander-Reiten quivers are coherent
(see Section \ref{sec:2} for definition).
Every coherent cycle-finite algebra $A$ is triangular, and hence the (geometric) Tits quadratic form $q_A$ and the (homological) Euler form
$\chi_A$ of $A$ are well defined. For a vector $\bf d$ in the Grothendieck group $K_0(A)$ of $A$ with nonnegative coordinates, we denote by
${\rm mod}_A({\bf d})$ the affine variety of $A$-modules of dimension vector $\bf d$ and by $G({\bf d})$ the corresponding product of general linear
groups acting on ${\rm mod}_A({\bf d})$ in such a way that the $G({\bf d})$-orbits in $\modu_A({\bf d})$ correspond to the isomorphism classes
of $A$-modules of dimension vector $\bf d$.
The third main aim of the article is to establish a common bound on the numbers of discrete indecomposable modules in each dimension vector over
cycle-finite algebras, generalizing results proved in \cite{[SZ]} for strongly simply connected algebras of polynomial growth.

For basic background from the representation theory of algebras we refer to the books \cite{[ASS]}, \cite{[ARS]}, \cite{[Ri]}, \cite{[SS1]},
\cite{[SS2]}, \cite{[SY]}.

\section{Preliminaries}\label{sec:2}

In this section we recall some concepts and results from the representation theory of algebras important for further considerations.

Let $A$ be an algebra (which by our assumption is basic and connected). Then there is an isomorphism $A\cong KQ/I$ of $K$-algebras, where $KQ$
is the path algebra of the {\it Gabriel quiver} $Q=Q_A$ of $A$ and $I$ is an admissible ideal of $KQ$. Equivalently, $A=KQ/I$ may be considered as
a $K$-category whose class of objects is the set $Q_0$ of vertices of $Q$, and the set of morphisms $A(x,y)$ from $x$ to $y$ is the quotient
of $K$-space $KQ(x,y)$, formed by the $K$-linear combinations of paths in $Q$ from $x$ to $y$, by the subspace $I(x,y)=KQ(x,y)\cap I$. We shall
identify an algebra $A$ with its $K$-category. Moreover, the module category $\modu A$ may be identified with the category $\rep_K(Q,I)$ of finite
dimensional $K$-linear representations of the bound quiver $(Q,I)$. An algebra $A$ with $Q_A$ acyclic (without oriented cycles) is said to be
{\it triangular}. A full subcategory $C$ of $A$ is said to be {\it convex} if any path in $Q_A$ with source and target in $Q_C$ lies entirely
in $Q_C$. Recall also that the {\it Jacobson radical} $\rad(\modu A)$ of the module category $\modu A$ is the ideal of $\modu A$ generated by all
noninvertible morphisms in $\ind A$. Then the {\it infinite radical} $\rad^{\infty}(\modu A)$ of $\modu A$ is the intersection of all powers
$\rad^i(\modu A)$, $i\geq 1$, of $\rad(\modu A$). A {\it path of length} $t\geq 1$ in $\modu A$ is a sequence of nonzero nonisomorphisms
\[ M_0 \buildrel {f_1}\over {\hbox to 6mm{\rightarrowfill}} M_1 \to \cdots \to M_{t-1}\buildrel {f_t}\over {\hbox to 6mm{\rightarrowfill}} M_t \]
and modules $M_0, M_1, \ldots, M_t$ in $\ind A$. Such a path is said to be {\it finite} if $f_1, \ldots, f_t$ do not belong to $\rad^{\infty}(\modu A)$,
and otherwise {\it infinite}. Moreover, if $M_0\cong M_t$ then the path is called a {\it cycle} of length $t$. A module $M$ from $\ind A$ is called
{\it directing} if it does not lie on a cycle in $\modu A$. For a module $M$ in $\modu A$, we denote by $\vdim M$ its {\it dimension vector}
$(\dim_KM(i))_{i\in Q_0}$. The {\it support} $\supp M$ of a module $M$ in $\modu A$ is the full subcategory of $A$ given by all vertices $i$ of $Q_A$
such that $M(i)\neq 0$. A module $M$ in $\modu A$ with $\supp M=A$ is said to be {\it sincere}. Recall also that the Grothendieck group
$K_0(A)=K_0(\modu A)$ is isomorphic to ${\Bbb Z}^{Q_0}$.

Let $A$ be an algebra and $K[x]$ the polynomial algebra in one variable $x$. Following \cite{[Dro]} $A$ is said to be {\it tame} if, for any dimension $d$,
there exists a finite number of $K[x]-A$-bimodules $M_i$, $1\leq i\leq n_d$, which are finitely generated and free as left $K[x]$-modules, and all but
a finite number of isoclasses of indecomposable $A$-modules of dimension $d$ are of the form $K[x]/(x-\lambda)\otimes_{K[x]}M_i$ for some $\lambda\in K$
and some $i\in\{1, \ldots ,n_d\}$. Let $\mu_A(d)$ be the least number of $K[x]-A$-bimodules $M_i$ satisfying the above condition for $d$. Then $A$
is said to be of {\it polynomial growth} (respectively, {\it domestic}) if there exists a positive integer $m$ such that $\mu_A(d)\leq d^m$
(respectively, $\mu_A(d)\leq m$) for any $d\geq 1$ (see \cite{[CB2]}, \cite{[Sk2]}). Recall that from the validity of the second Brauer-Thrall conjecture
we know that $A$ is representation-finite if and only if $\mu_A(d)=0$ for any $d\geq 1$.

The {\it Tits form} of a triangular algebra $A=KQ/I$ is the integral quadratic form $q_A : {\Bbb Z}^{Q_0}\to {\Bbb Z}$, defined, for
${\bf x}=(x_i)_{i\in Q_0}\in {\Bbb Z}^{Q_0}$, by
\[ q_{A}({\bf x}) = \sum_{i \in Q_0}x_{i}^{2} - \sum_{(i \rightarrow j)\in Q_{1}}x_{i}x_{j} + \sum_{i,j \in Q_{0}}r(i,j)x_{i}x_{j}, \]
where $Q_1$ is the set of arrows in $Q$ and $r(i,j)$ is the cardinality of $L\cap I(i,j)$, for a minimal set of generators
$L\subset\bigcup_{i,j\in Q_0}I(i,j)$ of the ideal $I$ (see \cite{[Bo2]}). Moreover, the {\it Euler form} of $A$ is the integral quadratic form
$\chi_A : {\Bbb Z}^{Q_0}\to {\Bbb Z}$ defined in \cite[2.4]{[Ri]} such that for any module $M$ in $\modu A$, we have
\[ \chi_A(\vdim M)=\sum_{i=0}^{\infty}(-1)^i\dim_K\Ext_A^i(M,M). \]
Observe that $A$ is of finite global dimension, because $A$ is triangular. It is also known that if $\gl\dim A\leq 2$ then $q_A=\chi_A$
(see \cite{[Bo2]}). Finally, it is known (see \cite{[Pe1]}) that, if $A$ is tame, then $q_A$ is {\it weakly nonnegative}, that is, $q_A({\bf x})\geq 0$
for all ${\bf x}\in{\Bbb N}^{Q_0}$. Unfortunately, the reverse implication is not true in general. However, it has been proved recently in \cite{[BPS]}
that a strongly simply connected algebra $A$ is tame if and only if the Tits form $q_A$ is weakly nonnegative. Recall also that a triangular algebra $A$
is called {\it strongly simply connected} \cite{[Sk3]} if the first Hochschild cohomology $H^1(C,C)$ of every convex subcategory $C$ of $A$ vanishes.

We need also special types of components of the Auslander-Reiten quivers of algebras.

Recall from \cite{[DR]}, \cite{[Ri]} that a translation quiver $\Gamma$ is called a {\it tube} if it contains a cyclical path and if its underlying
topological space is homeomorphic to $S^{1} \times {\Bbb R}^{+}$, where $S^1$ is the unit circle and ${\Bbb R}^+$ is the nonnegative real line.
A tube has only two types of arrows: arrows pointing to infinity and arrows pointing to the mouth. Tubes containing neither projective vertices
nor injective vertices are called {\it stable}, and are as follows. For the infinite quiver
\[ \mathbb{A}_{\infty}: 0\rightarrow 1\rightarrow 2\rightarrow \cdots\;\;\,\, \]
the translation quiver $\mathbb{Z}\mathbb{A}_{\infty}$ is of the form

\vspace{0.3cm}
\unitlength1cm
\begin{center}

\begin{picture}(7.7,2.7)
\put(0.4,2.75){\scriptsize($i-1,0$)}
\put(2.4,2.75){\scriptsize($i,0$)}
\put(4.4,2.75){\scriptsize($i+1,0$)}
\put(6.4,2.75){\scriptsize($i+2,0$)}
\put(1.4,1.75){\scriptsize($i-1,1$)}
\put(3.4,1.75){\scriptsize($i,1$)}
\put(5.4,1.75){\scriptsize($i+1,1$)}
\put(2.4,0.75){\scriptsize($i-1,2$)}
\put(4.4,0.75){\scriptsize($i,2$)}

\multiput(1,2.6)(1,-1){3}{\vector(1,-1){0.6}}
\multiput(3,2.6)(1,-1){3}{\vector(1,-1){0.6}}
\multiput(5,2.6)(1,-1){2}{\vector(1,-1){0.6}}
\multiput(7,2.6)(1,-1){1}{\vector(1,-1){0.6}}

\multiput(0,2)(1,1){1}{\vector(1,1){0.6}}
\multiput(1,1)(1,1){2}{\vector(1,1){0.6}}
\multiput(2,0)(1,1){3}{\vector(1,1){0.6}}
\multiput(4,0)(1,1){3}{\vector(1,1){0.6}}

\multiput(-0.2,1.8)(-0.15,-0.15){3}{$\cdot$}
\multiput(0.8,0.8)(-0.15,-0.15){3}{$\cdot$}
\multiput(1.8,-0.2)(-0.15,-0.15){3}{$\cdot$}

\multiput(7.7,1.7)(0.15,-0.15){3}{$\cdot$}
\multiput(6.7,0.7)(0.15,-0.15){3}{$\cdot$}
\multiput(5.7,-0.3)(0.15,-0.15){3}{$\cdot$}

\multiput(3.8,-0.3)(0,-0.15){3}{$\cdot$}

\end{picture}
\end{center}

\vspace{0.4cm}

\noindent with the translation $\tau$ given by $\tau(i,j) = (i-1,j)$ for $i \in \Bbb Z$, $j \in \Bbb N$. For each $r \geq 1$, denote by
${\Bbb Z\Bbb A}_{\infty}/(\tau^r)$ the translation quiver $\Gamma$ obtained from ${\Bbb Z\Bbb A}_{\infty}$ by identifying each vertex $(i,j)$ of
${\Bbb Z\Bbb A}_{\infty}$ with the vertex $\tau^r(i,j)$ and each arrow $x \rightarrow y$ in ${\Bbb Z\Bbb A}_{\infty}$ with the arrow
$\tau^r x \rightarrow \tau^r y$. The translation quiver ${\Bbb Z}{\Bbb A}_{\infty}/(\tau^{r})$ is called the {\it stable tube of rank $r$}.
The $\tau$-orbit of a stable tube $\Gamma$ formed by all vertices having exactly one immediate predecessor (equivalently, successor) is called the
{\it mouth} of $\Gamma$. A tube $\Gamma$ without injective vertices (respectively, without projective vertices) is called a {\it ray tube}
(respectively, {\it coray tube}).

Let $A$ be an algebra. A component $\mathcal C$ of $\Gamma_A$ is called {\it regular} if $\mathcal C$ contains neither a projective module nor an
injective module, and {\it semiregular} if $\mathcal C$ does not contain both a projective and an injective module. It has been shown in \cite{[Li1]}
and \cite{[Zh]} that a regular component $\mathcal C$ of $\Gamma_A$ contains an oriented cycle if and only if $\mathcal C$ is a stable tube.
Moreover, Liu proved in \cite{[Li2]} that a semiregular component $\mathcal C$ of $\Gamma_A$ contains an oriented cycle if and only if $\mathcal C$
is a ray or coray tube. A component $\mathcal P$ of $\Gamma_A$ is called {\it postprojective} if $\mathcal P$ is acyclic and every module in $\mathcal P$
lies in the $\tau_A$-orbit of a projective module. Dually, a component $\mathcal Q$ of $\Gamma_A$ is called {\it preinjective} if $\mathcal Q$ is acyclic
and every module in $\mathcal Q$ lies in the $\tau_A$-orbit of an injective module.
A component $\Gamma$ of $\Gamma_A$ is said to be {\it coherent} if the following two conditions are satisfied:

(C1) For each projective module $P$ in $\Gamma$ there is an infinite sectional path $P~=~X_1 \rightarrow X_2 \rightarrow \cdots \rightarrow X_i
\rightarrow X_{i+1} \rightarrow X_{i+2} \rightarrow \cdots$ (that is, $X_{i} \neq \tau_A X_{i+2}$ for any $i \geq 1$) in $\Gamma$.

(C2) For each injective module $I$ in $\Gamma$ there is  an infinite sectional path $\cdots \rightarrow Y_{j+2} \rightarrow Y_{j+1} \rightarrow
Y_j \rightarrow \cdots \rightarrow Y_2 \rightarrow Y_1 = I$ (that is, $Y_{j+2} \neq \tau_A Y_{j}$ for any $j \geq 1$) in $\Gamma$.

\noindent Further, a component $\Gamma$ of $\Gamma_A$ is said to be {\it almost cyclic} if all but finitely many modules of $\Gamma$ lie on oriented
cycles in $\Gamma_A$, so contained entirely in $\Gamma$. We note that the stable tubes, ray tubes and coray tubes of $\Gamma_A$ are special types of
almost cyclic coherent components. In general, it has been proved in \cite{[MS1]} that a component $\Gamma$ of $\Gamma_A$ is almost cyclic and coherent
if and only if $\Gamma$ is a {\it generalized multicoil}, which can be obtained from a family of stable tubes by a sequence of admissible operations
(see Section 4 for algebras having such components). A component $\Gamma$ of $\Gamma_A$ is said to be {\it almost acyclic} if all but finitely many
modules of $\Gamma$ are acyclic (do not lie on oriented cycles in $\Gamma_A$, hence in $\Gamma$), and {\it acyclic} if all modules of $\Gamma$ are
acyclic. Finally, following \cite{[Sk5]} a~component $\mathcal C$ of $\Gamma_A$ is said to be {\it generalized standard} if
$\rad_A^{\infty}(X,Y)=0$ for all modules $X$ and $Y$ from $\mathcal C$. It has been proved in \cite{[Sk5]} that every generalized standard component
$\mathcal C$ of $\Gamma_A$ is {\it almost periodic}, that is, all but finitely many $\tau_A$-orbits in $\mathcal C$ are periodic. Clearly, the
postprojective and preinjective components are acyclic, and the Auslander-Reiten quivers of representation-finite algebras are almost acyclic.
Moreover, these components are generalized standard (see \cite{[Sk6]}). General results on almost acyclic components and related algebras have
been proved by Reiten and third named author in \cite{[RS1]}, \cite{[RS2]} (see Section \ref{sec:5}).
For a component $\mathcal C$ of $\Gamma_A$, we denote by $\ann_A({\mathcal C})$ the annihilator of $\mathcal C$ in $A$, that is, the intersection of the
annihilators $\ann_A(X)=\{a\in A \mid Xa=0\}$ of all modules $X$ in $\mathcal C$. We note that $\mathcal C$ is a component of $\Gamma_{A/\ann_A({\mathcal C})}$.
Moreover, if $\ann_A({\mathcal C})=0$, $\mathcal C$ is said to be a {\it faithful component} of $\Gamma_A$.
By the {\it support} of a subquiver $\Gamma$ of $\Gamma_A$ we mean the full subcategory $\supp\Gamma$ of $A$ given by the supports $\supp M$ of all modules
$M$ in $\Gamma$, and, if $\supp\Gamma = A$ then $\Gamma$ is said to be {\it sincere}. We note that a faithful component $\mathcal C$ of $\Gamma_A$ is sincere.

\section{Semiregular components of cycle-finite algebras} \label{sec:3}

In this section we recall the shapes of the Auslander-Reiten quivers of representation-infinite tilted algebras of Euclidean type and tubular algebras,
as well as results from \cite{[Sk8]} on semiregular components of the Auslander-Reiten quivers of cycle-finite algebras, important for further
considerations.

By a {\it tame concealed algebra} we mean a tilted algebra $C=\End_H(T)$, where $H$ is the path algebra $K\Delta$ of a quiver $\Delta$ of
Euclidean type (the underlying graph $\overline\Delta$ of $\Delta$ of type $\widetilde{\Bbb A}_m$ ($m\geq 1$), $\widetilde{\Bbb D}_n$ ($n\geq 4$),
or $\widetilde{\Bbb E}_p$ ($6\leq p \leq 8$)) and $T$ is a (multiplicity-free) postprojective tilting $H$-modules. The tame concealed algebras
have been described by quivers and relations by Bongartz \cite{[Bo3]} and Happel-Vossieck \cite{[HV]}. Recall also that the Auslander-Reiten quiver
$\Gamma_C$ of a tame concealed algebra $C$ is of the form
\[ \Gamma_C=\mathcal{P}^C\cup\mathcal{T}^C\cup\mathcal{Q}^C, \]
where $\mathcal{P}^C$ is a postprojective component containing all indecomposable projective $C$-modules, $\mathcal{Q}^C$ is a preinjective component
containing all indecomposable injective $C$-modules, and $\mathcal{T}^C$ is a $\mathbb{P}_1(K)$-family $\mathcal{T}^C_{\lambda}$,
${\lambda\in\mathbb{P}_1(K)}$, of pairwise orthogonal generalized standard stable tubes, all but a finite number of them of rank one. The ordering
from the left to right indicates that there are nonzero homomorphisms only from any of these classes to itself and to the classes to its right. We
refer to \cite[Chapter 4]{[Ri]} and \cite{[SS1]} for more details on the module categories of tame concealed algebras.

Recall also that, if $B$ is a representation-infinite tilted algebra of Euclidean type $\Delta$, then one of the following holds:

(a) $B$ is a domestic tubular extension of a tame concealed algebra $C$ and
\[ \Gamma_B=\mathcal{P}^B\cup\mathcal{T}^B\cup\mathcal{Q}^B, \]
where $\mathcal{P}^B=\mathcal{P}^C$ is the postprojective component of $\Gamma_C$, $\mathcal{T}^B$ is a $\mathbb{P}_1(K)$-family $\mathcal{T}^B_{\lambda}$,
${\lambda\in\mathbb{P}_1(K)}$, of pairwise orthogonal generalized standard ray tubes, obtained from the $\mathbb{P}_1(K)$-family $\mathcal{T}^C$
of stable tubes of $\Gamma_C$ by ray insertions, and $\mathcal{Q}^B$ is a preinjective component containing all indecomposable injective $B$-modules
and a section of type $\Delta$;

(b) $B$ is a domestic tubular coextension of a tame concealed algebra $C$ and
\[ \Gamma_B=\mathcal{P}^B\cup\mathcal{T}^B\cup\mathcal{Q}^B, \]
where $\mathcal{P}^B$ is the postprojective component containing all indecomposable projective $B$-modules and a section of type $\Delta$,
$\mathcal{T}^B$ is a $\mathbb{P}_1(K)$-family $\mathcal{T}^B_{\lambda}$, ${\lambda\in\mathbb{P}_1(K)}$, of pairwise orthogonal generalized standard
coray tubes, obtained from the $\mathbb{P}_1(K)$-family $\mathcal{T}^C$ of stable tubes of $\Gamma_C$ by coray insertions, and
$\mathcal{Q}^B=\mathcal{Q}^C$ is the preinjective component of $\Gamma_C$.
We refer to \cite[Chapter 4]{[Ri]} and \cite[Chapters XV-XVII]{[SS2]} for more details on the module categories of representation-infinite tilted
algebras of Euclidean type.

By a {\it tubular algebra} we mean a tubular extension (equivalently tubular coextension) of a tame concealed algebra of tubular type $(2,2,2,2)$,
$(3,3,3)$, $(2,4,4)$, or $(2,3,6)$, as defined by Ringel in \cite[Chapter 5]{[Ri]}. Recall that a tubular algebra $B$ admits two different tame
concealed convex subcategories $C_0$ and $C_{\infty}$ such that $B$ is a tubular extension of $C_0$, and a tubular coextension of $C_{\infty}$,
and the Auslander-Reiten quiver $\Gamma_B$ of $B$ is of the form
\[ \Gamma_B=\mathcal{P}_0^B\cup\mathcal{T}_0^B\cup\left(\bigcup_{q\in{\Bbb Q}^+}\mathcal{T}_q^B\right)\cup\mathcal{T}_{\infty}^B\cup\mathcal{Q}_{\infty}^B, \]
where $\mathcal{P}_0^B=\mathcal{P}^{C_0}$ is the postprojective component of $\Gamma_{C_0}$, $\mathcal{T}_0^B$ is a $\mathbb{P}_1(K)$-family
of pairwise orthogonal generalized standard ray tubes, obtained from the $\mathbb{P}_1(K)$-family $\mathcal{T}^{C_0}$ of stable tubes of $\Gamma_{C_0}$
by ray insertions, $\mathcal{Q}_{\infty}^B=\mathcal{Q}^{C_{\infty}}$ is the preinjective component of $\Gamma_{C_{\infty}}$,
$\mathcal{T}_{\infty}^B$ is a $\mathbb{P}_1(K)$-family of pairwise orthogonal generalized standard coray tubes, obtained from the
$\mathbb{P}_1(K)$-family $\mathcal{T}^{C_{\infty}}$ of stable tubes of $\Gamma_{C_{\infty}}$ by coray insertions, and, for each $q\in{\Bbb Q}^+$
(the set of positive rational numbers) $\mathcal{T}_q^B$ is a $\mathbb{P}_1(K)$-family of pairwise orthogonal generalized standard stable tubes.
We refer to \cite[Chapter 5]{[Ri]} for more details on the module categories of tubular algebras.

The following characterization of tame concealed and tubular algebras has been established in \cite[Theorem 4.1]{[Sk8]}.
\begin{theorem} \label{thm31}
Let $A$ be an algebra. The following statements are equivalent:
\begin{enumerate}
\renewcommand{\labelenumi}{\rm(\roman{enumi})}
\item $A$ is cycle-finite and $\Gamma_A$ admits a sincere stable tube.
\item $A$ is either tame concealed or tubular.
\end{enumerate}
\end{theorem}
Moreover, we have also the following consequence \cite[Theorem 4.3]{[Sk8]} of the above theorem, the tameness of cycle-finite algebras
\cite[Proposition 1.4]{[AS4]}, and a result of Crawley-Boevey \cite[Corollary E]{[CB1]}.
\begin{theorem} \label{thm32}
Let $A$ be a cycle-finite algebra. Then $A$ is of polynomial growth.
\end{theorem}
An algebra $A$ is said to be {\it minimal representation-infinite} if $A$ is representation-infinite but every proper convex subcategory of $A$
is representation-finite. Then we have the following characterization of minimal representation-infinite cycle-finite algebras established in
\cite[Corollary 4.4]{[Sk8]}.
\begin{theorem} \label{thm33}
Let $A$ be an algebra. The following statements are equivalent:
\begin{enumerate}
\renewcommand{\labelenumi}{\rm(\roman{enumi})}
\item $A$ is minimal representation-infinite and cycle-finite.
\item $A$ is a tame concealed algebra.
\end{enumerate}
\end{theorem}
We also have the following characterization of domestic cycle-finite algebras established in \cite[Theorem 5.1]{[Sk8]}.
\begin{theorem} \label{thm34}
Let $A$ be a cycle-finite algebra. The following statements are equivalent:
\begin{enumerate}
\renewcommand{\labelenumi}{\rm(\roman{enumi})}
\item $A$ is domestic.
\item $A$ does not contain a tubular algebra as a convex subcategory.
\item All but finitely many components of $\Gamma_A$ are stable tubes of rank one.
\end{enumerate}
\end{theorem}
In general, we have the following information on the Auslander-Reiten quiver of a cycle-finite algebra, which is a consequence of results proved in
\cite{[Sk6]}, \cite{[Sk8]}, and results on the shapes of regular components from \cite{[Li1]}, \cite{[Zh]}.
\begin{theorem} \label{thm35}
Let $A$ be a cycle-finite algebra. Then every regular component of $\Gamma_A$ is a generalized standard stable tube. In particular, all but finitely
many components of $\Gamma_A$ are stable tubes.
\end{theorem}
Finally, we give a complete description of semiregular components of the Auslander-Reiten quivers of cycle-finite algebras, established in
\cite[Proposition 3.3]{[Sk8]}.
\begin{theorem} \label{thm36}
Let $A$ be a cycle-finite algebra and $\mathcal C$ be a semiregular component of $\Gamma_A$ containing a projective module. Then $B=\supp{\mathcal C}$
is a convex subcategory of $A$ and one of the following holds:
\begin{enumerate}
\renewcommand{\labelenumi}{\rm(\roman{enumi})}
\item $B$ is a domestic tubular coextension of a tame concealed algebra and $\mathcal C$ is the postprojective component of $\Gamma_B$.
\item $B$ is either a domestic tubular extension of a tame concealed algebra or a tubular algebra, and $\mathcal C$ is a generalized standard ray
tube of $\Gamma_B$.
\end{enumerate}
\end{theorem}
\begin{theorem} \label{thm37}
Let $A$ be a cycle-finite algebra and $\mathcal C$ be a semiregular component of $\Gamma_A$ containing an injective module. Then $B=\supp{\mathcal C}$
is a convex subcategory of $A$ and one of the following holds:
\begin{enumerate}
\renewcommand{\labelenumi}{\rm(\roman{enumi})}
\item $B$ is a domestic tubular extension of a tame concealed algebra and $\mathcal C$ is the preinjective component of $\Gamma_B$.
\item $B$ is either a domestic tubular coextension of a tame concealed algebra or a tubular algebra, and $\mathcal C$ is a generalized standard coray
tube of $\Gamma_B$.
\end{enumerate}
\end{theorem}

\section{Tame generalized multicoil algebras} \label{sec:4}

In this section we introduce and exhibit basic properties of the class of tame generalized multicoil algebras, plying a prominent role in the
description of infinite cyclic components of the Auslander-Reiten quivers of cycle-finite algebras. This is the class of tame algebras among the
class of all algebras having a separating family of almost cyclic coherent components investigated in \cite{[MS2]}, \cite{[MS3]}.
Recall that a family ${\mathcal C}$ = $({\mathcal C}_{i})_{i \in I}$ of components of the Auslander-Reiten quiver $\Gamma_A$ of an algebra $A$
is called {\it separating} in $\modu A$ if the modules in $\ind A$ split into three disjoint classes ${\mathcal P}^A$, ${\mathcal C}^A={\mathcal C}$
and ${\mathcal Q}^A$ such that:

(S1) ${\mathcal C}^A$ is a sincere generalized standard family of components;

(S2) $\Hom_{A}({\mathcal Q}^A,{\mathcal P}^A) = 0$, $\Hom_{A}({\mathcal Q}^A,{\mathcal C}^A)=0$, $\Hom_{A}({\mathcal C}^A,{\mathcal P}^A) = 0$;

(S3) any morphism from ${\mathcal P}^A$ to ${\mathcal Q}^A$ factors through the additive category

$\add {\mathcal C}^A$ of ${\mathcal C}^A$. \\
We then say that ${\mathcal C}^A$ {\it separates} ${\mathcal P}^A$ from ${\mathcal Q}^A$ and write
$\Gamma_A$=${\mathcal P}^A \cup {\mathcal C}^A \cup {\mathcal Q}^A$. We note that then ${\mathcal P}^A$ and ${\mathcal Q}^A$ are uniquely determined
by ${\mathcal C}^A$ (see \cite[(2.1)]{[AST]} or \cite[(3.1)]{[Ri]}).

We also recall a characterization of generalized standard stable tubes of an Auslander-Reiten quiver established in \cite[Corollary 5.3]{[Sk5]}
(see also \cite[Lemma 3.1]{[Sk7]}.
\begin{proposition}\label{prop41}
Let $A$ be an algebra and $\Gamma$ be a stable tube of $\Gamma_A$. The following statements are equivalent:
\begin{enumerate}
\renewcommand{\labelenumi}{\rm(\roman{enumi})}
\item $\Gamma$ is generalized standard.
\item The mouth of $\Gamma$ consists of pairwise orthogonal bricks.
\item $\rad_A^{\infty}(X,X)=0$ for any module $X$ in $\Gamma$.
\end{enumerate}
\end{proposition}
\noindent Recall that a module $X$ in $\modu A$ is called a {\it brick} if $\End_A(X)\cong K$.

It has been proved in \cite[Theorem A]{[MS1]} that a connected component $\Gamma$ of an Auslander-Reiten quiver $\Gamma_A$ is almost cyclic and coherent
if and only if $\Gamma$ is a generalized multicoil, obtained from a family of stable tubes by a sequence of operations called admissible.
We recall the letter and simultaneously define the corresponding enlargements of algebras.

We start with the one-point extensions and one-point coextensions of algebras. Let $A$ be an algebra and $M$ be a module in $\modu A$.
Then the {\it one-point extension} of $A$ by $M$ is the matrix algebra
\[ A[M]=\left[\begin{matrix}A&0\\M & K\end{matrix}\right] = \left\{\left[\begin{matrix}a&0\\m&\lambda
\end{matrix}\right];\,\,\lambda\in K,\,\, a\in A,\,\,m\in M\right\} \]
with the usual addition and multiplication. The quiver $Q_{A[M]}$ of $A[M]$ contains the quiver $Q_A$ of $A$ as a convex subquiver, and there is an
additional (extension) vertex which is a source.
The $A[M]$-modules are usually identified with the triples $(V,X,\varphi)$, where $V$ is a $K$-vector space, $X$ an $A$-module and
$\varphi : V \rightarrow $ Hom$_{A}(M,X)$ is a $K$-linear map. An $A[M]$-linear map $(V,X,\varphi) \rightarrow (W,Y,\psi)$ is then identified
with a pair $(f,g)$, where $f : V \rightarrow W$ is $K$-linear, $g : X \rightarrow Y$ is $A$-linear and $\psi f = $ Hom$_{A}(M,g)\varphi$.
Dually, one defines also the {\it one-point coextension} of $A$ by $M$ as the matrix algebra
\[ [M]A=\left[\begin{matrix}K & 0\\D(M) & A \end{matrix}\right]. \]
For $K$ and $r\geq 1$, we denote by $T_r(K)$ the $r\times r$-lower triangular matrix algebra
\[ \left [\begin{matrix} K & 0 & 0 & \ldots & 0 & 0 \\ K & K & 0 & \ldots &
0 & 0 \\ K & K & K & \ldots & 0 & 0 \\ \vdots & \vdots & \vdots &
\ddots & \vdots & \vdots \\ K & K & K & \ldots & K & 0 \\ K & K & K & \ldots & K & K
\end{matrix} \right ] \]

Given a generalized standard component $\Gamma$ of $\Gamma_A$, and an indecomposable
module $X$ in $\Gamma$, the {\it support} ${\mathcal S}(X)$ of the functor
$\Hom_{A}(X,-)\!\mid_{\Gamma}$ is the $R$-linear category  defined as follows
\cite{[AS6]}.
Let ${\mathcal H}_{X}$ denote the full subcategory of $\Gamma$ consisting of
the indecomposable modules $M$ in $\Gamma$ such that $\Hom_{A}(X,M) \neq 0$,
and ${\mathcal I}_{X}$ denote the ideal of ${\mathcal H}_X$ consisting of the
morphisms $f : M \rightarrow N$ (with $M, N$ in ${\mathcal H}_X$) such that
$\Hom_{A}(X,f) = 0$. We define ${\mathcal S}(X)$ to be the quotient category
${\mathcal H}_X/{\mathcal I}_X$. Following the above convention, we usually identify
the $R$-linear category ${\mathcal S}(X)$ with its quiver.

From now on let $A$ be an algebra and $\Gamma$ be a family of generalized standard infinite components of $\Gamma_A$. For
an indecomposable brick $X$ in $\Gamma$, called the {\it pivot}, one defines five admissible
operations (ad~1)-(ad~5) and their dual (ad~1$^*$)-(ad~5$^*$) modifying the translation quiver
$\Gamma =(\Gamma,\tau)$ to a new translation quiver $(\Gamma',\tau')$ and the algebra $A$ to a new algebra $A'$,
depending on the shape of the support ${\mathcal S}(X)$ (see \cite[Section 2]{[MS1]} for the figures illustrating the
modified translation quivers $\Gamma'$).

(ad~1) Assume ${\mathcal S}(X)$ consists of an infinite sectional path
starting at $X$:
\[ X = X_0 \rightarrow X_1 \rightarrow X_2 \rightarrow \cdots \]
\noindent In this case, we let
$t \geq 1$ be a positive integer, $D=T_t(K)$ and $Y_1$, $Y_2$, $\ldots$, $Y_t$ denote the
indecomposable injective $D$-modules with $Y=Y_1$ the unique indecomposable
projective-injective $D$-module. We define the {\it modified algebra} $A'$
of $A$ to be the one-point extension
\[ A' = (A\times D)[X\oplus Y] \]
\noindent and the {\it modified translation quiver} $\Gamma'$ of $\Gamma$
to be obtained by inserting in $\Gamma$ the rectangle consisting of
the modules $Z_{ij} = \left(K,X_i\oplus Y_j,\left[\begin{matrix}
1\\1\end{matrix}\right]\right)$ for $i \geq 0$, $1 \leq j \leq t$, and
$X'_{i} = (K,X_i,1)$ for $i \geq 0$. The translation $\tau'$ of
$\Gamma'$ is defined as follows: $\tau'Z_{ij} = Z_{i-1,j-1}$ if $i
\geq 1, j \geq 2, \tau'Z_{i1} = X_{i-1}$ if $i \geq 1, \tau'Z_{0j}
= Y_{j-1}$ if $j \geq 2, Z_{01}$ is projective, $\tau'X'_{0} =
Y_{t}, \tau'X'_{i} = Z_{i-1,t}$ if $i \geq 1,
\tau'(\tau^{-1}X_{i}) = X'_{i}$ provided $X_{i}$ is not an
injective $A$-module, otherwise $X'_{i}$ is injective in
$\Gamma'$. For the remaining vertices of $\Gamma'$, $\tau'$
coincides with the translation of $\Gamma$, or $\Gamma_D$,
respectively.

If $t = 0$ we define the modified algebra $A'$ to be the one-point extension
$A' = A[X]$ and the modified translation quiver $\Gamma'$ to be the translation quiver obtained
from $\Gamma$ by inserting only the sectional path consisting of the vertices
$X'_{i}$, $i \geq 0$.

The non-negative integer $t$ is such that the number of infinite sectional
paths parallel to $X_0 \to X_1 \to X_2 \to \cdots $ in the inserted
rectangle equals $t + 1$. We call $t$ the {\it parameter} of the operation.

Since $\Gamma$ is a generalized standard family of components of $\Gamma_A$, we then have

\begin{lemma} \label{ad1}
$\Gamma'$ is a generalized standard family of components of $\Gamma_{A'}$.
\end{lemma}

In case $\Gamma$ is a stable tube, it is clear that any module on the mouth
of $\Gamma$ satisfies the condition for being a pivot for the above
operation. Actually, the above operation is, in this case, the tube
insertion as considered in \cite{[DR]}.

(ad~2) Suppose that ${\mathcal S}(X)$ admits two sectional paths starting at $X$,
one infinite and the other finite with at least one arrow:
$$Y_t \leftarrow \cdots \leftarrow Y_2 \leftarrow Y_1 \leftarrow X = X_0
\rightarrow X_1 \rightarrow X_2 \rightarrow \cdots $$ \noindent
where $t\geq 1$. In particular, $X$ is necessarily injective. We
define the {\it modified algebra} $A'$ of $A$ to be the one-point
extension $A' = A[X]$ and the {\it modified translation quiver} $\Gamma'$
of $\Gamma$ to be obtained by inserting in $\Gamma$ the rectangle
consisting of the modules $Z_{ij} = \left(K,X_i\oplus
Y_j,\left[\begin{matrix} 1\\1\end{matrix}\right]\right)$ for $i \geq 1$,
$1 \leq j \leq t$, and $X'_{i} = (K,X_i,1)$ for $i \geq 1$. The
translation $\tau'$ of $\Gamma'$ is defined as follows: $X'_0$ is
projective-injective, $\tau'Z_{ij} = Z_{i-1,j-1}$ if $i \geq 2, j
\geq 2, \tau'Z_{i1} = X_{i-1}$ if $i \geq 1, \tau'Z_{1j} =
Y_{j-1}$ if $j \geq 2, \tau'X'_{i} = Z_{i-1,t}$ if $i \geq 2,
\tau'X'_1 = Y_t, \tau'(\tau^{-1}X_{i}) = X'_{i}$ provided $X_{i}$
is not an injective $A$-module, otherwise $X'_{i}$ is injective in
$\Gamma'$. For the remaining vertices of $\Gamma', \tau'$
coincides with the translation $\tau$ of $\Gamma$.

The integer $t \geq 1$ is such that the number of infinite sectional
paths parallel to $X_0 \to X_1 \to X_2 \to \cdots$ in the inserted
rectangle equals $t + 1$. We call $t$ the {\it parameter} of the operation.

Since $\Gamma$ is a generalized standard family of components of $\Gamma_A$, we then have

\begin{lemma} \label{ad2}
$\Gamma'$ is a generalized standard family of components of $\Gamma_{A'}$.
\end{lemma}

(ad~3) Assume ${\mathcal S}(X)$ is the mesh-category of two parallel sectional
paths:
\[ \begin{matrix} Y_1 & \to & Y_2 & \to & \cdots & \to & Y_t \\
\uparrow & & \uparrow & & & & \uparrow \\
X=X_0 & \to & X_1 & \to & \cdots & \to & X_{t-1} & \to & X_t & \to\cdots
\end{matrix} \]
\noindent where $t \geq 2$. In particular, $X_{t-1}$ is necessarily injective.
Moreover, we consider the translation quiver $\overline{\Gamma}$ of $\Gamma$ obtained by deleting the arrows
$Y_{i}\rightarrow \tau_{A}^{-1}Y_{i-1}$. We assume that the union $\widehat{\Gamma}$ of connected components of
$\overline{\Gamma}$ containing the vertices $\tau_{A}^{-1}Y_{i-1}$, $2\leq i\leq t$, is a finite translation quiver.
Then $\overline{\Gamma}$
is a disjoint union of $\widehat{\Gamma}$ and a cofinite full translation subquiver $\Gamma^*$, containing the pivot $X$.
We define the {\it modified algebra} $A'$ of $A$ to be the one-point extension $A' = A[X]$ and the {\it
modified translation quiver} $\Gamma'$ of $\Gamma$ to be obtained from $\Gamma^*$ by
inserting the rectangle consisting of the modules
$Z_{ij} = \left(K,X_i\oplus Y_j,\left[\begin{matrix}
1\\1\end{matrix}\right]\right)$ for $i \geq 1$, $1 \leq j \leq t$, and
$X'_{i} = (K,X_i,1)$ for $i \geq 1$. The translation $\tau'$ of
$\Gamma'$ is defined as follows: $X'_0$ is projective,
$\tau'Z_{ij} = Z_{i-1,j-1}$ if $i \geq 2$, $2 \leq j \leq t$,
$\tau'Z_{i1} = X_{i-1}$ if $i \geq 1, \tau'X'_i = Y_i$ if $1 \leq
i \leq t$, $\tau'X'_{i} = Z_{i-1,t}$ if $i \geq t + 1$, $\tau'Y_j
= X'_{j-2}$ if $2 \leq j \leq t$, $\tau'(\tau^{-1}X_{i}) =
X'_{i}$, if $i \geq t$ provided $X_{i}$ is not injective in
$\Gamma$, otherwise $X'_{i}$ is injective in $\Gamma'$. For the
remaining vertices of $\Gamma', \tau'$ coincides with the
translation $\tau$ of $\Gamma^*$. We note that $X'_{t-1}$ is
injective.

The integer $t \geq 2$ is such that the number of infinite sectional
paths parallel to $X_0 \to X_1 \to X_2 \to \cdots$ in the inserted
rectangle equals $t + 1$. We call $t$ the {\it parameter} of the operation.

Since $\Gamma$ is a generalized standard family of components of $\Gamma_A$, we then have

\begin{lemma} \label{ad3}
$\Gamma'$ is a generalized standard family of components of $\Gamma_{A'}$.
\end{lemma}

(ad~4) Suppose that ${\mathcal S}(X)$ consists an infinite sectional path, starting at $X$
\[ X = X_0 \rightarrow X_1 \rightarrow X_2 \rightarrow \cdots \]
\noindent and
\[ Y=Y_1 \rightarrow Y_2 \rightarrow \cdots \rightarrow Y_t \]
\noindent with $t \geq 1$, be a finite sectional path in $\Gamma_A$. Let $r$ be a positive integer.
Moreover, we consider the translation quiver $\overline{\Gamma}$ of $\Gamma$ obtained by deleting the arrows
$Y_{i}\rightarrow \tau_{A}^{-1}Y_{i-1}$. We assume that the union $\widehat{\Gamma}$ of connected components of
$\overline{\Gamma}$ containing the vertices $\tau_{A}^{-1}Y_{i-1}$, $2\leq i\leq t$, is a finite translation quiver.
Then $\overline{\Gamma}$
is a disjoint union of $\widehat{\Gamma}$ and a cofinite full translation subquiver $\Gamma^*$, containing the pivot $X$.
For $r = 0$ we define the {\it modified algebra} $A'$ of $A$ to be the
one-point extension $A' = A[X\oplus Y]$ and the {\it modified
translation quiver} $\Gamma'$ of $\Gamma$ to be obtained from $\Gamma^*$ by inserting the rectangle consisting of
the modules $Z_{ij} = \left(K,X_i\oplus Y_j,\left[\begin{matrix} 1\\1\end{matrix}\right]\right)$
for $i \geq 0$, $1 \leq j \leq t$, and $X'_{i} = (K,X_i,1)$ for $i
\geq 1$. The translation $\tau'$ of $\Gamma'$ is defined as
follows: $\tau'Z_{ij} = Z_{i-1,j-1}$ if $i \geq 1, j \geq 2,
\tau'Z_{i1} = X_{i-1}$ if $i \geq 1, \tau'Z_{0j} = Y_{j-1}$ if $j
\geq 2, Z_{01}$ is projective, $\tau'X'_{0} = Y_{t}, \tau'X'_{i} =
Z_{i-1,t}$ if $i \geq 1, \tau'(\tau^{-1}X_{i}) = X'_{i}$ provided
$X_{i}$ is not injective in $\Gamma$, otherwise $X'_{i}$ is
injective in $\Gamma'$. For the remaining vertices of $\Gamma',
\tau'$ coincides with the translation of $\Gamma^*$.

For $r \geq 1$, let $G=T_r(K)$, $U_{1,t+1}$, $U_{2,t+1}$, $\ldots$,
$U_{r,t+1}$ denote the indecomposable projective $G$-modules,
$U_{r,t+1}$, $U_{r,t+2}$, $\ldots$, $U_{r,t+r}$ denote the
indecomposable injective $G$-modules, with $U_{r,t+1}$ the unique
indecomposable projective-injective $G$-module. We define the {\it
modified algebra} $A'$ of $A$  to be the triangular matrix algebra
of the form:
\[ A'=\left [\begin{matrix} A & 0 & 0 & \ldots & 0 & 0 \\ Y & K & 0 & \ldots &
0 & 0 \\ Y & K & K & \ldots & 0 & 0 \\ \vdots & \vdots & \vdots &
\ddots & \vdots & \vdots \\ Y & K & K & \ldots & K & 0 \\ X\oplus Y & K & K & \ldots & K & K
\end{matrix} \right ] \]
\noindent with $r+2$ columns and rows and the {\it modified
translation quiver} $\Gamma'$ of $\Gamma$ to be obtained from $\Gamma^*$ by inserting the rectangles consisting of
the modules $U_{kl} = Y_l\oplus U_{k,t+k}$ for $1\leq k\leq r$, $1\leq l\leq t$, and $Z_{ij} =
\left(K,X_i\oplus U_{rj},\left[\begin{matrix} 1\\1\end{matrix}\right]\right)$
for $i \geq 0$, $1 \leq j \leq t+r$, and $X'_{i} = (K,X_i,1)$ for
$i \geq 0$. The translation $\tau'$ of $\Gamma'$ is defined as
follows: $\tau'Z_{ij} = Z_{i-1,j-1}$ if $i \geq 1, j \geq 2,
\tau'Z_{i1} = X_{i-1}$ if $i \geq 1, \tau'Z_{0j} = U_{r,j-1}$ if
$2\leq j \leq t+r$, $Z_{01}, U_{k1}, 1 \leq k \leq r$ are
projective, $\tau'U_{kl} = U_{k-1,l-1}$ if $2 \leq k \leq r$, $2
\leq l \leq t+r$, $\tau'U_{1l} = Y_{l-1}$ if $2 \leq l \leq t+1$,
$\tau'X'_{0} = U_{r,t+r}, \tau'X'_{i} = Z_{i-1,t+r}$ if $i \geq 1,
\tau'(\tau^{-1}X_{i}) = X'_{i}$ provided $X_{i}$ is not injective
in $\Gamma$, otherwise $X'_{i}$ is injective in $\Gamma'$. For the
remaining vertices of $\Gamma', \tau'$ coincides with the
translation of $\Gamma^*$, or $\Gamma_{G}$, respectively.

We note that the quiver $Q_{A'}$ of
$A'$ is obtained from the quiver of the double one-point extension $A[X][Y]$
by adding a path of length $r+1$ with source at the extension vertex of
$A[X]$ and sink at the extension vertex of $A[Y]$.

The integers $t \geq 1$ and $r\geq 0$ are such that the number of infinite
sectional paths parallel to $X_0 \to X_1 \to X_2 \to \cdots$ in the inserted
rectangles equals $t+r+1$. We call $t+r$ the {\it parameter} of the operation.

Since $\Gamma$ is a generalized standard family of components of $\Gamma_A$, we then have

\begin{lemma} \label{ad4}
$\Gamma'$ is a generalized standard family of components of $\Gamma_{A'}$.
\end{lemma}

(fad~1) Assume ${\mathcal S}(X)$ consists of a finite sectional path starting
at $X$:
\[ X = X_0 \rightarrow X_1 \rightarrow X_2 \rightarrow \cdots \rightarrow X_s \]
\noindent where $s \geq 0$ and $X_s$ is injective. Let $t \geq 1$ be a
positive integer, $D=T_t(K)$ and $Y_1$, $Y_2$, $\ldots$,
$Y_t$ denote the indecomposable injective $D$-modules with $Y=Y_1$ the
unique indecomposable projective-injective $D$-module. We define the
{\it modified algebra} $A'$ of $A$ to be the one-point extension
\[ A' = (A\times D)[X\oplus Y] \]
and the {\it modified translation quiver} $\Gamma'$ of $\Gamma$ to be
obtained by inserting in $\Gamma$ the rectangle consisting of the
modules $Z_{ij} = \left(K,X_i\oplus Y_j,\left[\begin{matrix}
1\\1\end{matrix}\right]\right)$ for $0\leq i \leq s$, $1 \leq j \leq t$,
$X'_{i} = (K,X_i,1)$ for $0\leq i \leq s$, $Y'_j = (K,Y_j,1)$ for
$1\leq j\leq t$, and $W=S_p$, where $p$ is the extension vertex of
$A[X]$. The translation $\tau'$ of $\Gamma'$ is defined as
follows: $\tau'Z_{ij} = Z_{i-1,j-1}$ if $i \geq 1, j \geq 2,
\tau'Z_{i1} = X_{i-1}$ if $i \geq 1, \tau'Z_{0j} = Y_{j-1}$ if $j
\geq 2, Z_{01}$ is projective, $\tau'X'_{0} = Y_{t}, \tau'X'_{i} =
Z_{i-1,t}$ if $i \geq 1, \tau'(\tau^{-1}X_{i}) = X'_{i}$ provided
$X_{i}$ is not injective in $\Gamma$, otherwise $X'_{i}$ is
injective in $\Gamma'$, $\tau'Y'_1 = X_s$, $\tau'Y'_j = Z_{s,j-1}$
if $2 \leq j \leq t$, $\tau'W = Z_{st}$. For the remaining
vertices of $\Gamma', \tau'$ coincides with the translation of
$\Gamma$, or $\Gamma_D$, respectively. If $t = 0$ we define the
modified algebra $A'$ to be the one-point extension $A' = A[X]$
and the modified translation quiver $\Gamma'$ to be the component obtained
from $\Gamma$ by inserting only the sectional path consisting of
the vertices $X'_{i}$, $0\leq i \leq s$, and $W$.

Observe that for $s = 0 = t$ the new translation
quiver $\Gamma'$ is obtained from $\Gamma$ by adding the projective-injective
vertex $X'_0$ and the vertex $W$ with $\tau'W = X_0$.

(fad~2) Suppose that ${\mathcal S}(X)$ admits two finite sectional paths starting
at $X$, each of them with at least one arrow:
\[ Y_t \leftarrow \cdots \leftarrow Y_2 \leftarrow Y_1 \leftarrow X = X_0 \rightarrow X_1 \rightarrow X_2 \rightarrow \cdots \rightarrow X_s \]
\noindent where $t\geq 1$ and $s\geq 1$. In particular, $X$
and $X_s$ are necessarily injective. We define the {\it modified
algebra} $A'$ of $A$ to be the one-point extension $A' = A[X]$ and
the {\it modified translation quiver} $\Gamma'$ of $\Gamma$ to be obtained
by inserting in $\Gamma$ the rectangle consisting of the modules
$Z_{ij} = \left(K,X_i\oplus Y_j,\left[\begin{matrix}
1\\1\end{matrix}\right]\right)$ for $1\leq i \leq s$, $1 \leq j \leq t$,
$X'_{i} = (K,X_i,1)$ for $1\leq i\leq s$, $Y'_j = (K,Y_j,1)$ for
$1\leq j\leq t$, and $W=S_p$, where $p$ is the extension vertex of
$A[X]$. The translation $\tau'$ of $\Gamma'$ is defined as
follows: $X'_0$ is projective-injective, $\tau'Z_{ij} =
Z_{i-1,j-1}$ if $i \geq 2, j \geq 2, \tau'Z_{i1} = X_{i-1}$ if $i
\geq 1, \tau'Z_{1j} = Y_{j-1}$ if $j \geq 2, \tau'X'_{i} =
Z_{i-1,t}$ if $i \geq 2, \tau'X'_1 = Y_t, \tau'(\tau^{-1}X_{i}) =
X'_{i}$ provided $X_{i}$ is not injective in $\Gamma$, otherwise
$X'_{i}$ is injective in $\Gamma'$, $\tau'Y'_1 = X_s$, $\tau'Y'_j
= Z_{s,j-1}$ if $2 \leq j \leq t$, $\tau'W = Z_{st}$. For the
remaining vertices of $\Gamma', \tau'$ coincides with the
translation $\tau$ of $\Gamma$.

(fad~3) Assume ${\mathcal S}(X)$ is the mesh-category of two finite parallel sectional
paths:
\[ \begin{matrix} Y_1 & \to & Y_2 & \to & \cdots & \to & Y_t \\ \uparrow & & \uparrow & & & & \uparrow \\
X=X_0 & \to & X_1 & \to & \cdots & \to & X_{t-1} & \to & X_t & \to\cdots\to X_s
\end{matrix} \]
\noindent where $s \geq t-1$, $t \geq 2$. In particular, $X_{t-1}$
and $X_s$ are necessarily injective. We define the {\it modified
algebra} $A'$ of $A$ to be the one-point extension $A' = A[X]$ and
the {\it modified translation quiver} $\Gamma'$ of $\Gamma$ to be obtained
by inserting in $\Gamma$ the rectangle consisting of the modules
$Z_{ij} = \left(K,X_i\oplus Y_j,\left[\begin{matrix}
1\\1\end{matrix}\right]\right)$ for $1\leq i \leq s$, $1 \leq j \leq t$,
$X'_{i} = (K,X_i,1)$ for $1\leq i \leq s$, $Y'_j = (K,Y_j,1)$ for
$1\leq j\leq t$, and $W=S_p$, where $p$ is the extension vertex of
$A[X]$. The translation $\tau'$ of $\Gamma'$ is defined as
follows: $X'_0$ is projective, $\tau'Z_{ij} = Z_{i-1,j-1}$ if $i
\geq 2, 2 \leq j \leq t, \tau'Z_{i1} = X_{i-1}$ if $i \geq 1,
\tau'X'_i = Y_i$ if $1 \leq i \leq t, \tau'X'_{i} = Z_{i-1,t}$ if
$i \geq t + 1, \tau'Y_j = X'_{j-2}$ if $2 \leq j \leq t$,
$\tau'(\tau^{-1}X_{i}) = X'_{i}$, if $i \geq t$ provided $X_{i}$
is not injective in $\Gamma$, otherwise $X'_{i}$ is injective in
$\Gamma'$. In both cases, $X'_{t-1}$ is injective, $\tau'Y'_1 =
X_s$, $\tau'Y'_j = Z_{s,j-1}$ if $2 \leq j \leq t$, $\tau'W =
Z_{st}$. For the remaining vertices of $\Gamma', \tau'$ coincides
with the translation $\tau$ of $\Gamma$. Observe that for $s = t -
1$ we have $Z_{tt} = Y'_t$ and $X'_t = W$.

(fad~4) Suppose that ${\mathcal S}(X)$ consists of a finite sectional path  starting at $X$:
\[ X = X_0 \rightarrow X_1 \rightarrow X_2 \rightarrow \cdots \rightarrow X_s \]
\noindent with $s \geq 1$ and $X_s$ injective, and
\[ Y=Y_1 \rightarrow Y_2 \rightarrow \cdots \rightarrow Y_t \]
\noindent $t \geq 1$,  be a finite sectional path in $\Gamma_A$.
Let $r$ be a positive integer. For $r = 0$ we define the
{\it modified algebra} $A'$ of $A$ to be the one-point
extension $A' = A[X\oplus Y]$ and the {\it modified translation quiver} $\Gamma'$ of
$\Gamma$ to be obtained by inserting in $\Gamma$ the rectangle consisting of
the modules $Z_{ij} = \left(K,X_i\oplus Y_j,\left[\begin{matrix} 1\\1\end{matrix}\right]\right)$
for $0\leq i \leq s$, $1\leq j\leq t$, $X'_{i} = (K,X_i,1)$ for
$0\leq i \leq s$, $Y'_j = (K,Y_j,1)$ for $1\leq j\leq t$, and $W=S_p$,
where $p$ is the extension vertex of $A[X]$.
The translation $\tau'$ of $\Gamma'$ is defined as follows: $\tau'Z_{ij} =
Z_{i-1,j-1}$ if $i \geq 1, j \geq 2, \tau'Z_{i1} = X_{i-1}$ if
$i \geq 1, \tau'Z_{0j} = Y_{j-1}$ if $j \geq 2, Z_{01}$ is
projective, $\tau'X'_{0} = Y_{t}, \tau'X'_{i} = Z_{i-1,t}$
if $i \geq 1, \tau'(\tau^{-1}X_{i}) = X'_{i}$ provided $X_{i}$
is not injective in $\Gamma$, otherwise $X'_{i}$ is injective in
$\Gamma'$, $\tau'Y'_1 = X_s$,
$\tau'Y'_j = Z_{s,j-1}$ if $2 \leq j \leq t$, $\tau'W = Z_{st}$.
For the remaining vertices of $\Gamma', \tau'$
coincides with the translation of $\Gamma$.

For $r \geq 1$, let $G=T_r(K)$, $U_{1,t+1}$, $U_{2,t+1}$, $\ldots$, $U_{r,t+1}$ denote the
indecomposable projective $G$-modules,
$U_{r,t+1}$, $U_{r,t+2}$, $\ldots$, $U_{r,t+r}$ denote the indecomposable injective
$G$-modules, with $U_{r,t+1}$ the unique indecomposable projective-injective
$G$-module. We define the {\it modified algebra} $A'$ of $A$  to be
the triangular matrix algebra of the form:
\[ A'=\left [\begin{matrix} A & 0 & 0 & \ldots & 0 & 0 \\ Y & K & 0 & \ldots &
0 & 0 \\ Y & K & K & \ldots & 0 & 0 \\ \vdots & \vdots & \vdots &
\ddots & \vdots & \vdots \\ Y & K & K & \ldots & K & 0 \\ X\oplus Y & K & K & \ldots & K & K
\end{matrix} \right ] \]
\noindent with $r+2$ columns and rows and the {\it modified translation quiver}
$\Gamma'$ of $\Gamma$ to be obtained by inserting in $\Gamma$ the rectangles
consisting of the modules
$U_{kl} = Y_l\oplus U_{k,t+k}$ for $1\leq k\leq r$, $1\leq l\leq t$,
$Z_{ij} = \left(K,X_i\oplus U_{rj},\left[\begin{matrix} 1\\1\end{matrix}\right]\right)$ for
$0\leq i \leq s$, $1 \leq j \leq t+r$, $X'_{i} = (K,X_i,1)$ for
$0\leq i \leq s$, $Y'_j = (K,U_{rj},1)$ for $1\leq j\leq t+r$, and $W=S_p$,
where $p$ is the extension vertex of $A[X]$.
The translation $\tau'$ of $\Gamma'$ is defined as follows:
$\tau'Z_{ij} = Z_{i-1,j-1}$ if $i \geq 1, j \geq 2,
\tau'Z_{i1} = X_{i-1}$ if $i \geq 1, \tau'Z_{0j} = U_{r,j-1}$ if
$2\leq j \leq t+r$, $Z_{01}, U_{k1}, 1 \leq k \leq r$ are
projective, $\tau'U_{kl} = U_{k-1,l-1}$ if $2 \leq k \leq r$,
$2 \leq l \leq t+r$, $\tau'U_{1l} = Y_{l-1}$ if $2 \leq l \leq t+1$,
$\tau'X'_{0} = U_{r,t+r}, \tau'X'_{i} = Z_{i-1,t+r}$
if $i \geq 1, \tau'(\tau^{-1}X_{i}) = X'_{i}$ provided $X_{i}$
is not injective in $\Gamma$, otherwise $X'_{i}$ is injective in
$\Gamma'$, $\tau'Y'_1 = X_s$,
$\tau'Y'_j = Z_{s,j-1}$ if $2 \leq j \leq t+r$, $\tau'W = Z_{s,t+r}$.
For the remaining vertices of $\Gamma', \tau'$
coincides with the translation of $\Gamma$, or $\Gamma_{G}$, respectively.

(ad~5) We define the {\it modified algebra} $A'$ of $A$ to be the iteration
of the extensions described in the definitions of the admissible operations
(ad~1), (ad~2), (ad~3), (ad~4), and their finite versions corresponding to
the operations (fad~1), (fad~2), (fad~3) and (fad~4). The
{\it modified translation quiver} $\Gamma'$ of $\Gamma$ is obtained in
the following three steps: first we are doing on $\Gamma$ one of
the operations (fad~1), (fad~2) or (fad~3), next a finite number (possibly
empty) of the operation (fad~4) and finally the operation (ad~4), and in such
a way that the sectional paths starting from all the new projective vertices
have a common cofinite (infinite) sectional subpath.

Since $\Gamma$ is a generalized standard family of components of $\Gamma_A$, we then have

\begin{lemma} \label{ad5}
$\Gamma'$ is a generalized standard family of components of $\Gamma_{A'}$.
\end{lemma}

Finally, together with each of the admissible operations (ad~1), (ad~2),
(ad~3), (ad~4) and (ad~5), we consider its dual, denoted by (ad~1$^{*}$),
(ad~2$^{*}$), (ad~3$^{*}$), (ad~4$^{*}$) and (ad~5$^*$). These ten operations are
called the {\it admissible operations}.
Following \cite{[MS1]} a connected translation quiver $\Gamma$ is said to be a
{\it generalized multicoil} if $\Gamma$ can be obtained
from a finite family ${\mathcal T}_1, {\mathcal T}_2, \ldots, {\mathcal T}_s$ of stable tubes
by an iterated application of admissible operations
(ad~1), (ad~1$^{*}$), (ad~2), (ad~2$^{*}$), (ad~3), (ad~3$^{*}$), (ad~4),
(ad~4$^{*}$), (ad~5) or (ad~5$^{*}$). If $s = 1$, such a translation quiver
$\Gamma$ is said to be a {\it generalized coil}.
The admissible operations of types (ad~1), (ad~2), (ad~3), (ad~1$^{*}$), (ad~2$^{*}$) and (ad~3$^{*}$) have been
introduced in \cite{[AS4]}, \cite{[AS6]}, \cite{[AST]}, and the admissible operations (ad~4) and (ad~4$^{*}$) for $r=0$ in \cite{[Ma1]}.

Observe that any stable tube is trivially a generalized coil. A
{\it tube} (in the sense of \cite{[DR]}) is a generalized coil
having the property that each admissible operation in the sequence
defining it is of the form (ad~1) or (ad~1$^*$). If we apply only
operations of type (ad~1) (respectively, of type (ad 1$^*$)) then
such a generalized coil is called a {\it ray tube} (respectively,
a {\it coray tube}). Observe that a generalized coil without
injective (respectively, projective) vertices is a ray tube
(respectively, a coray tube). A {\it quasi-tube} (in the sense of
\cite{[Sk2]}) is a generalized coil having the property that each
of the admissible operations in the sequence defining it is of
type (ad~1), (ad~1$^{*}$), (ad~2) or (ad~2$^{*}$). Finally,
following \cite{[AS6]} a coil is a generalized coil having the
property that each of the admissible operations in the sequence
defining it is one of the forms (ad~1), (ad~1$^{*}$), (ad~2),
(ad~2$^{*}$), (ad~3) or (ad~3$^{*}$). We note that any generalized
multicoil $\Gamma$ is a coherent translation quiver with trivial
valuations and its cyclic part $_c\Gamma$ (the translation
subquiver of $\Gamma$ obtained by removing from $\Gamma$ all
acyclic vertices and the arrows attached to them) is infinite,
connected and cofinite in $\Gamma$, and so $\Gamma$ is almost
cyclic.

Let $C$ be the product $C_1\times\ldots\times C_m$ of a family $C_1, \ldots, C_m$ of tame concealed algebras and ${\mathcal T}^C$ the disjoin union
${\mathcal T}^{C_1}\cup\ldots\cup{\mathcal T}^{C_m}$ of $\mathbb{P}_1(K)$-families ${\mathcal T}^{C_1}, \ldots, {\mathcal T}^{C_m}$ of pairwise
orthogonal generalized standard stable tubes of $\Gamma_{C_1}, \ldots, \Gamma_{C_m}$, respectively. Following \cite{[MS2]}, we say that an algebra
$A$ is a {\it generalized multicoil enlargement} of $C_1, \ldots, C_m$ if $A$ is obtained from $C$ by an iteration of admissible operations
of types (ad~1)-(ad~5) and (ad~1$^{*}$)-(ad~5$^{*}$) performed either on stable tubes of ${\mathcal T}^C$ or on generalized multicoils obtained from
stable tubes of ${\mathcal T}^C$ by means of operations done so far. It follows from \cite[Corollary B]{[MS2]} that then $A$ is a triangular algebra,
and hence the Tits and Euler forms $q_A$ and $\chi_A$ of $A$ are defined. In fact, in \cite{[MS2]} generalized multicoil enlargements of finite
families of arbitrary concealed canonical algebras have been introduced and investigated. But in the tame case we may restrict to the generalized
multicoil enlargements of tame concealed algebras. Namely, we have the following consequence of \cite[Theorems A and F]{[MS2]}.
\begin{theorem} \label{thm42}
Let $A$ be an algebra. The following statements are equivalent:
\begin{enumerate}
\renewcommand{\labelenumi}{\rm(\roman{enumi})}
\item $A$ is tame and $\Gamma_A$ admits a separating family of almost cyclic coherent components.
\item $A$ is a tame generalized multicoil enlargement of a finite family of tame concealed algebras.
\item $A$ is a generalized multicoil enlargement of a finite family of tame concealed algebras and the Tits form $q_A$ is weakly nonnegative.
\end{enumerate}
\end{theorem}
From now on, by a {\it tame generalized multicoil algebra} we mean a connected tame generalized multicoil enlargement of a finite family of tame
concealed algebras. The following consequence of \cite[Theorems C and F]{[MS2]} describes the structure of the Auslander-Reiten quivers of tame
generalized multicoil algebras.
\begin{theorem} \label{thm43}
Let $A$ be a tame generalized multicoil algebra obtained from a family $C_1, \ldots, C_m$ of tame concealed algebras. There are convex subcategories
$A^{(l)}=A_1^{(l)}\times\ldots\times A_m^{(l)}$ and $A^{(r)}=A_1^{(r)}\times\ldots\times A_m^{(r)}$ of $A$ such that the following statement hold:
\begin{enumerate}
\renewcommand{\labelenumi}{\rm(\roman{enumi})}
\item For each $i\in\{1, \ldots, m\}$, $A_i^{(l)}$ is a domestic tubular or tubular coextension of the tame concealed algebra $C_i$.
\item For each $i\in\{1, \ldots, m\}$, $A_i^{(r)}$ is a domestic tubular or tubular extension of the tame concealed algebra $C_i$.
\item The Auslander-Reiten quiver $\Gamma_A$ of $A$ is of the form
$$\Gamma_A={\mathcal P}^A \cup {\mathcal C}^A \cup {\mathcal Q}^A,$$
where ${\mathcal C}^A$ is a family of generalized multicoils separating ${\mathcal P}^A$ from ${\mathcal Q}^A$ such that:
\end{enumerate}
\begin{enumerate}
\renewcommand{\labelenumi}{\rm(\arabic{enumi})}
\item ${\mathcal C}^A$ is obtained from the $\mathbb{P}_1(K)$-families ${\mathcal T}^{C_1}, \ldots, {\mathcal T}^{C_m}$ of stable tubes of
$\Gamma_{C_1}, \ldots, \Gamma_{C_m}$ by admissible operations corresponding to the admissible operations leading from $C_1, \ldots, C_m$ to $A$;
\item ${\mathcal P}^A$ is the disjoint union ${\mathcal P}^{A_1^{(l)}}\cup\ldots\cup{\mathcal P}^{A_m^{(l)}}$, where, for each $i\in\{1, \ldots, m\}$,
${\mathcal P}^{A_i^{(l)}}$ is either the postprojective component of $\Gamma_{A_i^{(l)}}$, if $A_i^{(l)}$ is tilted of Euclidean type, or
${\mathcal P}^{A_i^{(l)}}={\mathcal P}_0^{A_i^{(l)}}\cup {\mathcal T}_0^{A_i^{(l)}}\cup\left(\bigcup_{q\in{\Bbb Q}^+}{\mathcal T}_q^{A_i^{(l)}}\right)$,
if $A_i^{(l)}$ is a tubular algebra;
\item ${\mathcal Q}^A$ is the disjoint union ${\mathcal Q}^{A_1^{(r)}}\cup\ldots\cup{\mathcal Q}^{A_m^{(r)}}$, where, for each $i\in\{1, \ldots, m\}$,
${\mathcal Q}^{A_i^{(r)}}$ is either the preinjective component of $\Gamma_{A_i^{(r)}}$, if $A_i^{(r)}$ is tilted of Euclidean type, or
${\mathcal Q}^{A_i^{(r)}}=\left(\bigcup_{q\in{\Bbb Q}^+}{\mathcal T}_q^{A_i^{(r)}}\right)\cup{\mathcal T}_{\infty}^{A_i^{(r)}}\cup {\mathcal Q}_{\infty}^{A_i^{(r)}}$,
if $A_i^{(r)}$ is a tubular algebra.
\end{enumerate}
\end{theorem}
In particular, we have the following consequence of Theorems \ref{thm31} and \ref{thm43}.
\begin{corollary}
Let $A$ be a tame generalized multicoil algebra. Then $A$ is cycle-finite.
\end{corollary}
Further, as a consequence of Theorems \ref{thm34} and \ref{thm43}, we obtain the following fact.
\begin{corollary}
Let $A$ be a tame generalized multicoil algebra and $\Gamma_A={\mathcal P}^A \cup {\mathcal C}^A \cup {\mathcal Q}^A$ the canonical decomposition
of $\Gamma_A$. The following statements are equivalent:
\begin{enumerate}
\renewcommand{\labelenumi}{\rm(\roman{enumi})}
\item $A$ is domestic.
\item ${\mathcal P}^A$ is a disjoint union of postprojective components of Euclidean type and ${\mathcal Q}^A$ is a disjoint union of
preinjective components of Euclidean type.
\end{enumerate}
\end{corollary}
Moreover, the following consequence of \cite[Theorem E]{[MS2]} describes the homological properties of modules over tame generalized multicoil algebras.
\begin{theorem}
Let $A$ be a tame generalized multicoil algebra and $\Gamma_A={\mathcal P}^A \cup {\mathcal C}^A \cup {\mathcal Q}^A$ the canonical decomposition
of $\Gamma_A$ described above. Then the following statements hold:
\begin{enumerate}
\renewcommand{\labelenumi}{\rm(\roman{enumi})}
\item $\pd_AX\leq 1$ for any module $X$ in ${\mathcal P}^A$.
\item $\id_AX\leq 1$ for any module $X$ in ${\mathcal Q}^A$.
\item $\pd_AX\leq 2$ and $\id_AX\leq 2$ for any module $X$ in ${\mathcal C}^A$.
\item $\gldim A\leq 3$.
\end{enumerate}
\end{theorem}

\section{Tame generalized double tilted algebras} \label{sec:5}

In this section we introduce and describe basic properties of the class of tame generalized double tilted algebras, which is the class of tame algebras
in the class of generalized double tilted algebras investigated in \cite{[RS1]}, \cite{[RS2]}, \cite{[Sk11]}.

Let $H$ be a hereditary algebra, $T$ a (multiplicity-free) tilting module in $\modu H$ and $B=\End_H(T)$ the associated tilted algebra.
Then $T$ induces the torsion pair $({\mathcal T}(T),{\mathcal F}(T))$ in $\modu H$, with the torsion class
${\mathcal T}(T)=\{M\in\modu H; \Ext^1_A(T,M)=0\}$ and the torsion-free class ${\mathcal T}(F)=\{N\in\modu H; \Hom_H(T,N)=0\}$, and the torsion pair
$({\mathcal X}(T),{\mathcal Y}(T))$ in $\modu B$, with the torsion class
${\mathcal X}(T)=\{X\in\modu B; X\otimes_BT=0\}$ and the torsion-free class ${\mathcal Y}(T)=\{Y\in\modu H; \Tor^B_1(Y,T)=0\}$. Then, by the
Brenner-Butler theorem, the functors $\Hom_A(T,-)$ and $-\otimes_BT$ induce quasi-inverse equivalence between ${\mathcal T}(T)$ and ${\mathcal Y}(T)$,
and the functors $\Ext^1_A(T,-)$ and $\Tor^B_1(-,T)$ induce quasi-inverse equivalence between ${\mathcal F}(T)$ and ${\mathcal X}(T)$ (see \cite{[BB]}, \cite{[HR]}).
Moreover, $({\mathcal X}(T),{\mathcal Y}(T))$ is a splitting torsion pair of $\modu B$, that is, every indecomposable module in $\modu B$ belongs
either to ${\mathcal X}(T)$ or ${\mathcal Y}(T)$. Further, the images $\Hom_H(T,I)$ of the indecomposable injective $H$-modules $I$ via the functor
$\Hom_H(T,I)$ form a section $\Sigma$ of an acyclic component ${\mathcal C}_T$ of $\Gamma_B$ such that $\Sigma$ is isomorphic to the opposite
quiver $Q_H^{\op}$ of the quiver $Q_H$ of $H$, any predecessor of $\Sigma$ in ${\mathcal C}_T$ lies in ${\mathcal Y}(T)$, and any proper successor
of $\Sigma$ in ${\mathcal C}_T$ lies in ${\mathcal X}(T)$. Therefore, the component ${\mathcal C}_T$ of $\Gamma_B$ connects the torsion-free part
${\mathcal Y}(T)$ with the torsion part ${\mathcal X}(T)$ along the section $\Sigma$, and hence ${\mathcal C}_T$ is called the {\it connecting component}
of $\Gamma_B$ determined by $T$.

The following theorem proved independently by Liu \cite{[Li3]} and Skowro\'nski \cite{[Sk4]} gives a handy criterion for an algebra to be a tilted
algebra.
\begin{theorem} \label{thm51}
An algebra $B$ is a tilted algebra if and only if $\Gamma_B$ contains a component $\mathcal C$ with a faithful section $\Sigma$ such that
$\Hom_B(U,\tau_BV)=0$ for all modules $U, V$ from $\Sigma$. Moreover, in this case, the direct sum $T$ of all modules on $\Sigma$ is a tilting
$B$-module, $H=\End_B(T)$ is a hereditary algebra, $T^*=D(_HT)$ is a tilting $H$-module with $B\cong\End_H(T^*)$, and $\mathcal C$ is the
connecting component ${\mathcal C}_{T^*}$ of $\Gamma_B$ determined by $T^*$.
\end{theorem}
The general shape of the Auslander-Reiten quiver of a tilted algebra has been described by Kerner in \cite{[Ke]}. We will describe only the
Auslander-Reiten quivers of tame tilted algebras, which are exactly the cycle-finite tilted algebras.
\begin{theorem} \label{thm52}
Let $H = K\Delta$ be a hereditary algebra, $T$ a tilting $H$-module, and assume that the associated tilted algebra $B = \End_H(T)$ is tame.
Then the connecting component ${\mathcal C}_T$ of $\Gamma_B$ determined by $T$ admits a finite (possibly empty) family of pairwise disjoint
translation subquivers ${\mathcal D}_1^{(l)}, \dots, {\mathcal D}_m^{(l)}$, ${\mathcal D}_1^{(r)}, \dots, {\mathcal D}_n^{(r)}$
such that the following statements hold.
\begin{enumerate}
\renewcommand{\labelenumi}{\rm(\roman{enumi})}
\item For each $i \in \{ 1, \dots, m \}$, there exists an isomorphism of translation quivers ${\mathcal D}_i^{(l)} \cong {\Bbb N}\Delta_i^{(l)}$, where
$\Delta_i^{(l)}$ is a connected convex subquiver  of $\Delta$ of Euclidean type and ${\mathcal D}_i^{(l)}$  is closed under predecessors in
${\mathcal C}_T$.
\item For each $j \in \{ 1, \dots, n \}$, there exists an isomorphism of translation quivers ${\mathcal D}_j^{(r)} \cong {(-\Bbb N)}\Delta_j^{(r)}$, where
$\Delta_j^{(r)}$ is a connected convex subquiver  of $\Delta$ of Euclidean type and ${\mathcal D}_j^{(r)}$  is closed under successors in
${\mathcal C}_T$.
\item  All but finitely many modules of ${\mathcal C}_T$ lie in ${\mathcal D}_1^{(l)} \cup \dots \cup {\mathcal D}_m^{(l)} \cup
  {\mathcal D}_1^{(r)} \cup \dots \cup {\mathcal D}_n^{(r)}$.
\item For each $i \in \{ 1, \dots, m \}$, there exists a tilted algebra $B_i^{(l)} = \End_{H_i^{(l)}} (T_i^{(l)})$, where $H_i^{(l)}$ is the path
algebra $K \Delta_i^{(l)}$, $T_i^{(l)}$ is a tilting $H_i^{(l)}$-module without nonzero preinjective direct summands, $B_i^{(l)}$ is a quotient
algebra of $B$, and ${\mathcal D}_i^{(l)}$ coincides with the torsion-free part ${\mathcal Y}(T_i^{(l)}) \cap {\mathcal C}_{T_i^{(l)}}$ of the connecting
component ${\mathcal C}_{T_i^{(l)}}$ of $\Gamma_{B_i^{(l)}}$ determined by $T_i^{(l)}$.
\item For each $j \in \{ 1, \dots, n \}$, there exists a tilted algebra $B_j^{(r)} = \End_{H_j^{(r)}} (T_j^{(r)})$, where $H_j^{(r)}$ is the path
algebra $K \Delta_j^{(r)}$, $T_j^{(r)}$ is a tilting $H_j^{(r)}$-module without nonzero postprojective direct summands, $B_j^{(r)}$ is a quotient
algebra of $B$, and ${\mathcal D}_j^{(r)}$ coincides with the torsion part ${\mathcal X}(T_j^{(r)}) \cap {\mathcal C}_{T_j^{(r)}}$ of the connecting
component ${\mathcal C}_{T_j^{(r)}}$ of $\Gamma_{B_j^{(r)}}$ determined by $T_j^{(r)}$.
\item ${\mathcal Y}(T) = \add ({\mathcal Y}(T_1^{(l)}) \cup \dots \cup {\mathcal Y}(T_m^{(l)})\cup ({\mathcal Y}(T) \cap {\mathcal C}_T))$.
\item ${\mathcal X}(T) = \add (({\mathcal X}(T) \cap {\mathcal C}_T) \cup {\mathcal X}(T_1^{(r)}) \cup \dots \cup {\mathcal X}(T_n^{(r)}))$.
\item The Auslander-Reiten quiver $\Gamma_B$ of $B$ has the disjoint union decomposition
$$\Gamma_B =\left( \bigcup_{i=1}^m {\mathcal Y} \Gamma_{B_i^{(l)}} \right) \cup {\mathcal C}_T \cup
\left( \bigcup_{j=1}^n {\mathcal X} \Gamma_{B_j^{(r)}}  \right),$$
where
\end{enumerate}
\begin{enumerate}
\renewcommand{\labelenumi}{\rm(\alph{enumi})}
\item For each $i \in \{ 1, \dots, m \}$, ${\mathcal Y}\Gamma_{B_i^{(l)}}$ is the union of all components of $\Gamma_{B_i^{(l)}}$
contained entirely in ${\mathcal Y}({T_i^{(l)}})$, and hence consists of a unique postprojective component ${\mathcal P}^{B_i^{(l)}}$
and a ${\Bbb P}_1(K)$-family ${\mathcal T}^{B_i^{(l)}} = ({\mathcal T}_{\lambda}^{B_i^{(l)}})_{\lambda \in {\Bbb P}_1 (K)}$
of pairwise orthogonal generalized standard ray tubes;
\item For each $j \in \{ 1, \dots, n \}$, ${\mathcal X}\Gamma_{B_j^{(r)}}$ is the union of all components of $\Gamma_{B_j^{(r)}}$
contained entirely in ${\mathcal X}({T_j^{(r)}})$, and hence consists of a unique preinjective component ${\mathcal Q}^{B_j^{(r)}}$
and a ${\Bbb P}_1(K)$-family ${\mathcal T}^{B_j^{(r)}} = ({\mathcal T}_{\lambda}^{B_j^{(r)}})_{\lambda \in {\Bbb P}_1 (K)}$
of pairwise orthogonal generalized standard coray tubes.
\end{enumerate}
\end{theorem}
The following theorem follows from \cite{[Bo2]}(part (i)) and \cite[p.376]{[Ri]}(parts (ii) and (iii)).
\begin{theorem}\label{thm53}
Let $A$ be a cycle-finite algebra, $X$ a directing module in $\modu A$, and $B=\supp X$. Then the following statements hold:
\begin{enumerate}
\renewcommand{\labelenumi}{\rm(\roman{enumi})}
\item $B$ is a convex subcategory of $A$.
\item $B$ is a tame tilted algebra.
\item $X$ belongs to a connecting component of $\Gamma_B$.
\end{enumerate}
\end{theorem}
We refer to \cite{[Bo1]} and \cite{[Dra]} (respectively, \cite{[Pe2]} and \cite{[Pe3]}) for a classification of representation-finite
(respectively, tame representation-infinite) tame tilted algebras with sincere directing modules.

The class of tilted algebras was extended in \cite{[RS1]} to the class of double tilted algebras, and next in \cite{[RS2]} to the class of
generalized double tilted algebras, containing the class of all algebras of finite type, by extending the concept of a section to the concept
of a multisection.

Following \cite{[RS2]}, a full connected subquiver $\Delta$ of a component $\mathcal C$ of the Auslander-Reiten quiver $\Gamma_A$ of an algebra
$A$ is said to be a {\it multisection} if the following conditions are satisfied:
\begin{enumerate}
\renewcommand{\labelenumi}{\rm(\roman{enumi})}
\item $\Delta$ is almost acyclic.
\item $\Delta$ is convex in $\mathcal C$.
\item For each $\tau_A$-orbit $\mathcal O$ in $\mathcal C$, we have $1 \leq | \Delta \cap {\mathcal O} | < \infty$.
\item $| \Delta \cap {\mathcal O} | = 1$ for all but finitely many $\tau_A$-orbits $\mathcal O$ in $\mathcal C$.
\item No proper full convex subquiver of $\Delta$ satisfies \rm{(i)}--\rm{(iv)}.
\end{enumerate}

It has been proved in \cite[Theorem 2.5]{[RS2]} that a component $\mathcal C$ of $\Gamma_A$ is almost acyclic if and only if $\mathcal C$ admits
a multisection $\Delta$. Moreover, for a multisection $\Delta$ of a component $\mathcal C$ of $\Gamma_A$, the following full subquivers of
$\mathcal C$ were defined in \cite{[RS2]}:
\begin{enumerate}
\renewcommand{\labelenumi}{\rm(\roman{enumi})}
\item $\! \Delta^{\prime}_{l} = \{ X\! \in \!\Delta; \text{there is a nonsectional path in $\mathcal C$ from $X$ to a projective module $P$}\},$
\item $\! \Delta^{\prime}_{r} = \{ X\! \in \!\Delta; \text{there is a nonsectional path in $\mathcal C$ from an injective module $I$ to $X$}\},$
\item $
  \Delta^{\prime\prime}_{l} =
  \{ X \in \Delta^{\prime}_{l};
     \tau_A^{-1} X \notin \Delta^{\prime}_{l}
  \} , \qquad
  \Delta^{\prime\prime}_{r} =
  \{ X \in \Delta^{\prime}_{r};
     \tau_A X \notin \Delta^{\prime}_{r}
  \} ,
$
\item $
  \Delta_l =
  (\Delta \setminus \Delta^{\prime}_r)
       \cup \tau_A \Delta^{\prime\prime}_r
    , \quad
  \Delta_c = \Delta^{\prime}_l \cap \Delta^{\prime}_r
    , \quad
  \Delta_r = (\Delta \setminus \Delta^{\prime}_l)
        \cup \tau_A^{-1} \Delta^{\prime\prime}_l
  .
$
\end{enumerate}
Then $\Delta_l$ is called the {\it left part} of $\Delta$, $\Delta_r$ the {\it right part} of $\Delta$, and $\Delta_c$ the {\it core} of $\Delta$.
\begin{lemma}\label{lem54}
Let $A$ be an algebra, $\mathcal C$ a component of $\Gamma_A$ and $\Delta$ a multisection of $\mathcal C$. The following statements hold
\begin{enumerate}
\renewcommand{\labelenumi}{\rm(\roman{enumi})}
\item Every cycle of $\mathcal C$ lies in $\Delta_c$.
\item $\Delta_c$ is finite.
\item Every indecomposable module $X$ in $\mathcal C$ is in $\Delta_c$, or a predecessor of $\Delta_l$ or a successor of $\Delta_r$ in $\mathcal C$.
\item $\Delta$ is faithful if and only if $\mathcal C$ is faithful.
\end{enumerate}
\end{lemma}
Moreover, in \cite{[RS2]} a numerical invariant $w(\Delta)\in{\Bbb N}\cup\{\infty\}$ of a multisection $\Delta$ of $\mathcal C$, called the
{\it width} of $\Delta$, was introduced such that $\mathcal C$ is acyclic if and only if $w(\Delta)<\infty$, and $w(\Delta)=1$ if and only if
$\Delta$ is a section.

The following facts proved in \cite[Proposition 2.11]{[RS2]} show that the core and the width of a multisection of an almost cyclic component
$\mathcal C$ are uniquely determined by $\mathcal C$.
\begin{proposition} \label{prop55}
Let $A$ be an algebra, $\mathcal C$ a component of $\Gamma_A$ and $\Delta, \Sigma$ multisections of $\mathcal C$. Then $\Delta_c=\Sigma_c$ and
$w(\Delta)=w(\Sigma)$.
\end{proposition}

Following \cite{[RS2]}, an algebra $B$ is said to be a {\it generalized double tilted algebra} if the following conditions are satisfied:
\begin{enumerate}
\renewcommand{\labelenumi}{\rm(\arabic{enumi})}
\item $\Gamma_B$ admits a component $\mathcal C$ with a faithful multisection $\Delta$.
\item There exists a tilted quotient algebra $B^{(l)}$ of $B$ (not necessarily connected) such that $\Delta_l$ is a disjoint union of sections
of the connecting components of the connected parts of $B^{(l)}$ and the category of all predecessors of $\Delta_l$ in $\ind B$ coincides with
the category of all predecessors of $\Delta_l$ in $\ind B^{(l)}$.
\item There exists a tilted quotient algebra $B^{(r)}$ of $B$ (not necessarily connected) such that $\Delta_r$ is a disjoint union of sections
of the connecting components of the connected parts of $B^{(r)}$, and the category of all successors of $\Delta_r$ in $\ind B$ coincides with
the category of all successors of $\Delta_r$ in $\ind B^{(r)}$.
\end{enumerate}
Then $B^{(l)}$ is called a {\it left tilted part} of $B$ and $B^{(r)}$ a {\it right tilted part} of $B$.

The following generalization of Theorem \ref{thm51}, proved in \cite[Theorem 3.1]{[RS2]}, gives a handy criterion for an algebra to be a generalized double
tilted algebra.
\begin{theorem}\label{thm55}
Let $B$ be an algebra. The following conditions are equivalent:
\begin{enumerate}
\renewcommand{\labelenumi}{\rm(\roman{enumi})}
\item $B$ is a generalized double tilted algebra.
\item  The quiver $\Gamma_B$ admits a component $\mathcal C$ with a faithful multisection $\Delta$ such that $\Hom_B(U, \tau_B V) = 0$, for all modules
 $U \in \Delta_r$ and $V \in \Delta_l$.
\item The quiver $\Gamma_B$ admits a faithful generalized standard almost acyclic component $\mathcal C$.
\end{enumerate}
\end{theorem}
In particular, we obtain the following characterization of tame generalized double tilted algebras.
\begin{theorem}\label{thm56}
Let $B$ be a generalized double tilted algebra, $\mathcal C$ a faithful generalized standard almost cyclic component of $\Gamma_B$,
and $\Delta$ a multisection of $\mathcal C$. The following conditions are equivalent:
\begin{enumerate}
\renewcommand{\labelenumi}{\rm(\roman{enumi})}
\item $B$ is tame.
\item $B$ is cycle-finite.
\item $\Delta_l$ and $\Delta_r$ are disjoint unions of Euclidean quivers.
\item The tilted algebras $B^{(l)}$ and $B^{(r)}$ are tame.
\item The Auslander-Reiten quiver $\Gamma_B$ of $B$ has disjoint union decomposition
$$\Gamma_B = {\mathcal Y}\Gamma_{B^{(l)}}\cup {\mathcal C}\cup {\mathcal X}\Gamma_{B^{(r)}},$$
where
\end{enumerate}
\begin{enumerate}
\renewcommand{\labelenumi}{\rm(\alph{enumi})}
\item ${\mathcal Y}\Gamma_{B^{(l)}}$ is the union of all connected components of $\Gamma_{B^{(l)}}$ contained entirely in the torsion-free part
${\mathcal Y}(B^{(l)})$, and ${\mathcal Y}\Gamma_{B^{(l)}}$ is a disjoint union of postprojective components of Euclidean type and
${\Bbb P}_1(K)$-families of pairwise orthogonal generalized standard ray tubes.
\item ${\mathcal X}\Gamma_{B^{(r)}}$ is the union of all connected components of $\Gamma_{B^{(r)}}$ contained entirely in the torsion part
${\mathcal X}(B^{(r)})$, and ${\mathcal X}\Gamma_{B^{(r)}}$ is a disjoint union of preinjective components of Euclidean type and
${\Bbb P}_1(K)$-families of pairwise orthogonal generalized standard coray tubes.
\end{enumerate}
\end{theorem}
We end this section with an example of a tame generalized double tilted algebra, illustrating the above considerations.
\begin{example}
Let $B=KQ/I$, where $Q$ is the quiver
$$\xymatrix@C=12pt@R=6pt{
1 \\
&&6\ar[ldd]\\
2\\
&5\ar[luuu]\ar[lu]\ar[ld]\ar[lddd]&&7\ar[luu]\ar[ll]&8\ar@<2pt>[l]\ar@<-2pt>[l]\\
3\\
\\
4\\
}$$
and $I$ is the ideal of $KQ$ generated by all paths of $Q$ of length 2 (see \cite[Example 4.3]{[RS2]}). Then $B$ is a tame generalized double tilted algebra
of global dimension 4 and $\Gamma_B$ admits a generalized standard component $\mathcal C$ of the form
$$\xymatrix@C=10pt@R=2pt{
&{\phantom{I}}\ar[rddd]&&I_1\ar[rddd]\\
&&&&&P_6\ar[rdd]&&I_6\ar[rdd]&&P_8\ar@<2pt>[rdd]\ar@<-2pt>[rdd]&&{\phantom{P}}&\cdots \\
&{\phantom{I}}\ar[rd]&&I_2\ar[rd]\\
\cdots&&\bullet\ar[ruuu]\ar[ru]\ar[rd]\ar[rddd]&&S_5\ar[ruu]\ar[rdd]&&R\ar[ruu]\ar[rdd]&&S_7\ar@<2pt>[ruu]\ar@<-2pt>[ruu]&&\bullet\ar@<2pt>[ruu]\ar@<-2pt>[ruu]\\
&{\phantom{I}}\ar[ru]&&I_3\ar[ru]\\
&&&&&P_7\ar[ruu]\ar[rdd]&&I_5\ar[rdd]\ar[ruu] \\
&{\phantom{I}}\ar[ruuu]&&I_4\ar[ruuu] \\
&&&&S_6\ar[ruu]\ar[ruu]&&P_7/S_6\ar[ruu]&&S_6 \\
}$$
with a faithful multisection $\Delta$ formed by the indecomposable injective modules $I_1, I_2, I_3, I_4, I_5, I_6$ (at the vertices 1, 2, 3, 4, 5, 6),
the indecomposable projective modules $P_6, P_7, P_8$ (at the vertices 6, 7, 8), the simple modules $S_5, S_6, S_7$ (at the vertices 5, 6, 7) and
the modules $P_7/S_6$, $R=\tau^{-1}_BS_5=\tau_BS_7$. Then the left part $\Delta_l$ of $\Delta$ consists of the modules $I_1, I_2, I_3, I_4, S_5, P_6$
and is a section of the preinjective connecting component of the tame tilted algebra $B^{(l)}$ being the convex subcategory of $B$ given by the vertices
1, 2, 3, 4, 5, 6. The right part $\Delta_r$ of $\Delta$ consists of the modules $I_6, S_7, P_8$ and is a section of the postprojective connecting
component of the tame tilted algebra $B^{(r)}$ being the convex subcategory of $B$ given by the vertices 6, 7, 8. Moreover, the core
$\Delta_c$ of $\Delta$ consists of the modules $S_6, P_7, P_7/S_6, R, I_5$ and is the cyclic part of the Auslander-Reiten quiver $\Gamma_{B^{(c)}}$ of
the representation-finite convex subcategory $B^{(c)} = \supp\Delta_c$ of $B$ given by the vertices 5, 6, 7. We also note that every module in
$\ind B$ belongs to one of its full subcategories $\ind B^{(l)}$, $\ind B^{(c)}$, or $\ind B^{(r)}$.
\end{example}

\section{Cyclic components of cycle-finite algebras} \label{sec:6}

Let $A$ be an algebra. We denote by $_c\Gamma_A$ the translation subquiver of $\Gamma_A$, called the {\it cyclic part} of $\Gamma_A$, obtained  by
removing from $\Gamma_A$ all acyclic modules and the arrows attached to them. The connected components of $_c\Gamma_A$ are said to be
{\it cyclic components} of $\Gamma_A$ (see \cite{[MS1]}). The following result from \cite[Proposition 5.1]{[MS1]} will be very useful.
\begin{proposition} \label{prop62}
Let $A$ be an algebra and $X$, $Y$ be two cyclic modules of $\Gamma_A$. Then $X$ and $Y$ belong to the same cyclic component of $\Gamma_A$ if
and only if there is an oriented cycle in $\Gamma_A$ passing through $X$ and $Y$.
\end{proposition}
Moreover, we have the following property of the support algebras of cyclic components of the Auslander-Reiten quivers of cycle-finite algebras
(see \cite{[MPS]}).
\begin{proposition} \label{prop63}
Let $A$ be a cycle-finite algebra, $\Gamma$ a cyclic component of $\Gamma_A$ and $B=\supp\Gamma$. Then $B$ is a convex subcategory of $A$.
\end{proposition}

Let $A$ be an algebra and $\mathcal C$ be a component of $\Gamma_A$. We denote by $_l\mathcal C$ the {\it left stable part} of $\mathcal C$, obtained
by removing from $\mathcal C$ the $\tau_A$-orbits containing projective modules, and by $_r\mathcal C$ the {\it right stable part} of $\mathcal C$,
obtained by removing from $\mathcal C$ the $\tau_A$-orbits containing injective modules. We note that if $\mathcal C$ is an infinite component
of $\Gamma_A$ then $_l\mathcal C$ or $_r\mathcal C$ is not empty.

The following theorem from \cite[Theorem 1]{[MPS]} describes the supports of infinite cyclic components of the Auslander-Reiten quivers
of cycle-finite algebras.
\begin{theorem}\label{thm64}
Let $A$ be a cycle-finite algebra and $\Gamma$ an infinite cyclic component of $\Gamma_A$. Then there exist infinite full translation subquivers
$\Gamma_1, \ldots, \Gamma_r$ of $\Gamma$ such that the following statements hold.
\begin{enumerate}
\renewcommand{\labelenumi}{\rm(\roman{enumi})}
\item For each $i\in\{1,\ldots,r\}$, $\Gamma_i$ is a cyclic coherent full translation subquiver of $\Gamma_A$.
\item For each $i\in\{1,\ldots,r\}$, $B^{(i)}=\supp\Gamma_i$ is a tame generalized multicoil algebra and a quotient algebra of $A$.
\item $\Gamma_1, \ldots, \Gamma_r$ are pairwise disjoint full translation subquivers of $\Gamma$ and
$\Gamma^{cc}=\Gamma_1\cup\ldots\cup\Gamma_r$ is a maximal cyclic coherent and cofinite full translation subquiver of $\Gamma$.
\item $B(\Gamma\setminus\Gamma^{cc})=A/\ann_A(\Gamma\setminus\Gamma^{cc})$ is of finite representation type.
\end{enumerate}
\end{theorem}

The following example illustrates the above theorem.
\begin{example}
Let $A=KQ/I$, where $Q$ is the quiver
$$\xymatrix@C=12pt@R=12pt{
&&&9\ar[rd]^{\eta}&&16\ar[ld]_{\psi}\ar[rd]_{l}&&19\ar[rd]^{i}\ar[ll]_{j}\\
&&10\ar[ru]^{\xi}\ar[rd]_{\mu}\ar[dd]_{\pi}&&7\ar[dddd]^{\rho}&&17\ar[d]_{m}&&20\ar[d]^{h}\\
&&&8\ar[ru]_{\nu}&&&18&&21\ar[dd]^{g}\ar[rd]^{f}\\
0&1\ar[l]_{\theta}&2\ar[l]_{\omega}\ar[dd]_{\lambda}&&&&&&&22\ar[ld]^{e}\\
&&&5\ar[ld]_{\beta}&&&&&15\ar[d]^{d}\\
&&3&&6\ar[ld]^{\sigma}\ar[lu]_{\alpha}\ar[rd]_{\varphi}&&12\ar[ld]^{a}\ar[rd]_{b}&&14\ar[ld]^{c}\\
&&&4\ar[lu]^{\gamma}&&11&&13\\
}$$
and $I$ is the ideal in the path algebra $KQ$ of $Q$ generated by the elements $\alpha\beta - \sigma\gamma$, $\xi\eta - \mu\nu$,
$\pi\lambda - \xi\eta\rho\alpha\beta$, $\rho\varphi$, $\psi\rho$, $jl$, $dc$, $ed$, $gd$, $hg$, $hf$, $ih$. Then $A$ is a cycle-finite and
$\Gamma_A$ admits a component $\mathcal C$ of the form
$$\includegraphics[scale=0.60]{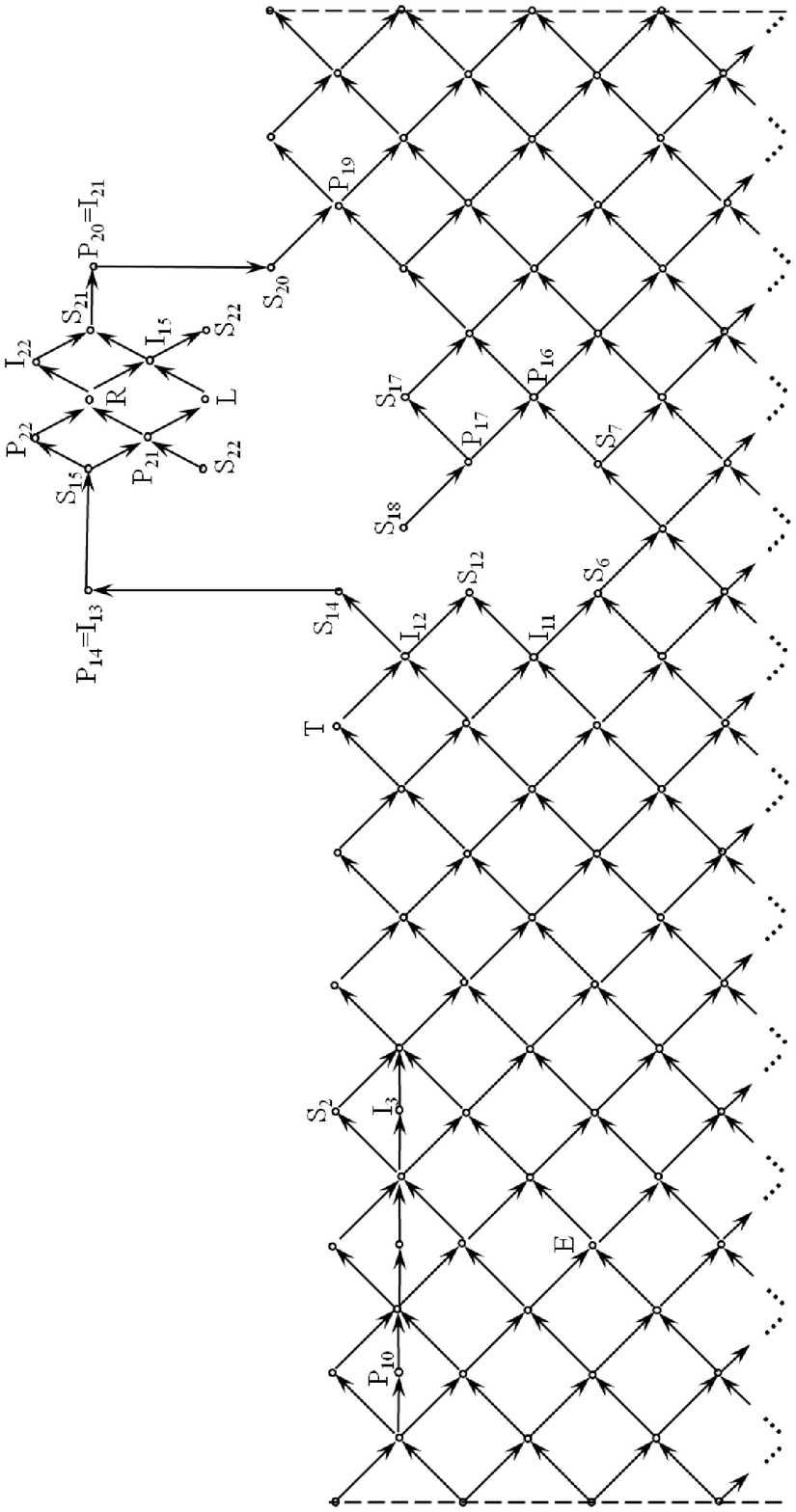}$$
The cyclic part $\Gamma$ of $\mathcal C$ is obtain by removing from $\mathcal C$ the (directing) modules $S_{12}, S_{17}, S_{18}, P_{17}$ and the
arrows attached to them. The maximal cyclic coherent part $\Gamma^{cc}$ of $\Gamma$ is the full translation subquiver of $\mathcal C$ obtained by
removing the modules $S_{12}$, $I_{12}$, $T$, $S_{14}$, $P_{14}=I_{13}$, $S_{15}$, $P_{21}$, $S_{22}$, $L$, $P_{22}$, $R$, $I_{15}$, $I_{22}$,
$S_{21}$, $P_{20}=I_{21}$, $S_{20}$, $S_{17}$, $P_{17}$, $S_{18}$ and the arrows attached to them. Further, $\Gamma^{cc}$ is the cyclic part
of the maximal almost cyclic coherent full translation subquiver $\Gamma^*$ of $\mathcal C$ obtained by removing the modules $P_{14}=I_{13}$,
$S_{15}$, $P_{21}$, $S_{22}$, $L$, $P_{22}$, $R$, $I_{15}$, $I_{22}$, $S_{21}$, $P_{20}=I_{21}$ and the arrows attached to them.

Let $B=A/\ann\Gamma$. Then $B=A/\ann\Gamma^*$, because $\ann\Gamma = \ann\Gamma^*$. Observe that $B=KQ_B/I_B$, where $Q_B$ is the full subquiver
of $Q$ given by all vertices of $Q$ except 15, 21, 22, and $I_B=I\cap KQ_B$. We claim that $B$ is a tame generalized multicoil algebra. Consider
the path algebra $C=K\Delta$ of the full subquiver of $Q$ given by the vertices 4, 5, 6, 7, 8, 9. Then $C$ is a hereditary algebra of Euclidean type
$\widetilde{\BD}_5$, and hence a tame concealed algebra. It is know that $\Gamma_C$ admits $\mathbb{P}_1(K)$-family $\mathcal{T}^{C}_{\lambda}$,
${\lambda\in\mathbb{P}_1(K)}$, of pairwise orthogonal generalized standard stable tubes, having a unique stable tube $\mathcal{T}_{\lambda}$ of rank
3 with the mouth formed by the modules $S_6=\tau_CS_7$, $S_7=\tau_CE$, $E=\tau_CS_6$, where $E$ is the unique indecomposable $C$-module with the
dimension vector
$\vdim E = {\begin{smallmatrix}1{\phantom{1}} \\ {\phantom{1}}1 \\ 1{\phantom{1}} \\ {\phantom{1}}{\phantom{1}} \\
1{\phantom{1}} \\ {\phantom{1}}1 \\ 1{\phantom{1}}
\end{smallmatrix}}$, (see \cite[Theorem XIII 2.9]{[SS1]}).

Then $B$ is the generalized multicoil enlargement of $C$, obtained by applications
\begin{itemize}
\item two admissible operations of types (ad~1$^*$) with the pivots $S_6$ and $S_{12}$, creating the vertices 11, 12, 13 and the arrows
$\varphi$, $a$, $b$, $c$;
\item two admissible operations of types (ad~1$^*$) with the pivots $E$ and $S_{2}$, creating the vertices 3, 2, 1, 0 and the arrows
$\beta$, $\gamma$, $\lambda$, $\omega$, $\theta$;
\item two admissible operations of types (ad~1) with the pivots $S_7$ and $S_{16}$, creating the vertices 16, 17, 18, 19, 20 and the arrows
$\psi$, $l$, $m$, $j$, $i$;
\item one admissible operation of type (ad~3) with the pivot the radical of $P_{10}$, creating the vertex 10 and the arrows $\xi$, $\mu$, $\pi$.
\end{itemize}
Then the left part $B^{(l)}$ of $B$ is the convex subcategory of $B$ (and of $A$) given by the vertices 0, 1, 2, 3, 4, 5, 6, 7, 8, 9, 11, 12, 13, 14,
and is a tilted algebra of Euclidean type $\widetilde{\BD}_{13}$ with the connecting postprojective component ${\mathcal P}^{B^{(l)}}$ containing all
indecomposable projective $B^{(l)}$-modules.
The right part $B^{(r)}$ of $B$ is the convex subcategory of $B$ (and of $A$) given by the vertices 4, 5, 6, 7, 8, 9, 16, 17, 18, 19, 20,
and is a tilted algebra of Euclidean type $\widetilde{\BD}_{10}$ with the connecting preinjective component ${\mathcal Q}^{B^{(r)}}$ containing all
indecomposable injective $B^{(r)}$-modules. We also note that the algebra $B(\Gamma\setminus\Gamma^{cc}) = A/\ann_A(\Gamma\setminus\Gamma^{cc})$ is
the representation-finite convex subcategory of $A$ given by the vertices 12, 13, 14, 15, 20, 21, 22. It follows from Theorem \ref{thm43} that the
Auslander-Reiten quiver $\Gamma_B$ of the generalized multicoil enlargement $B$ of $C$ is of the form
$$\Gamma_B={\mathcal P}^B \cup {\mathcal C}^B \cup {\mathcal Q}^B,$$
where ${\mathcal P}^B = {\mathcal P}^{B^{(l)}}$, ${\mathcal Q}^B = {\mathcal Q}^{B^{(r)}}$, and ${\mathcal C}^B$ is the $\mathbb{P}_1(K)$-family
${\mathcal C}^{B}_{\lambda}$, ${\lambda\in\mathbb{P}_1(K)}$, of pairwise orthogonal generalized multicoils such that ${\mathcal C}^B_1 = \Gamma^*$
and ${\mathcal C}^{B}_{\lambda}={\mathcal T}^{C}_{\lambda}$ for all ${\lambda\in\mathbb{P}_1(K)}\setminus\{1\}$. Hence $\Gamma_A$ is of the form
$$\Gamma_A={\mathcal P}^A \cup {\mathcal C}^A \cup {\mathcal Q}^A,$$
where ${\mathcal P}^A = {\mathcal P}^{B^{(l)}}$, ${\mathcal Q}^A = {\mathcal Q}^{B^{(r)}}$, and ${\mathcal C}^A$ is the $\mathbb{P}_1(K)$-family
${\mathcal C}^{A}_{\lambda}$, ${\lambda\in\mathbb{P}_1(K)}$, of pairwise orthogonal generalized standard components such that
${\mathcal C}^A_1 = {\mathcal C}$, ${\mathcal C}^{A}_{\lambda}={\mathcal T}^{C}_{\lambda}$ for all ${\lambda\in\mathbb{P}_1(K)}\setminus\{1\}$.
Moreover, we have
$$\Hom_A({\mathcal C}^{A},{\mathcal P}^{A})=0, \Hom_A({\mathcal Q}^{A},{\mathcal C}^{A})=0, \Hom_A({\mathcal Q}^{A},{\mathcal P}^{A})=0.$$
In particular, $A$ is a cycle-finite algebra with $(\rad^{\infty}(\modu A))^3=0$.
\end{example}

The following theorem from \cite[Theorem 2]{[MPS]} describes the supports of finite cyclic components of the Auslander-Reiten quivers
of cycle-finite algebras.
\begin{theorem} \label{thm65}
Let $A$ be a cycle-finite algebra and $\Gamma$ a finite cyclic component of $\Gamma_A$. Then the following statements hold.
\begin{enumerate}
\renewcommand{\labelenumi}{\rm(\roman{enumi})}
\item $B=\supp\Gamma$ is a tame generalized double tilted algebra.
\item $\Gamma$ is the core of the connecting component ${\mathcal C}_B$ of $\Gamma_B$.
\end{enumerate}
\end{theorem}

We note that if the core of an almost cyclic component of an Auslander-Reiten quiver $\Gamma_A$ is not empty that it contains a projective module
and an injective module. Then we obtain the following consequence of Theorem \ref{thm65}.
\begin{corollary} \label{cor66}
Let $A$ be a cycle-finite algebra. Then the number of finite cyclic components of $\Gamma_A$ is bounded by the rank of $K_0(A)$.
\end{corollary}

Observe also that for a cycle-finite algebra $A$ of infinite representation type there are infinitely many (infinite) cyclic components of $\Gamma_A$,
since $\Gamma_A$ contains infinitely many stable tubes (see Theorems \ref{thm31} and \ref{thm33}).

\section{The structure theorems} \label{sec:7}

Let $A$ be a cycle-finite algebra and $X$ a module in $\ind A$. Then $X$ is a directing module if and only if $X$ is an acyclic vertex
of $\Gamma_A$. Hence $X$ is nondirecting if and only if $X$ belongs to a cyclic component of $\Gamma_A$. Then the following structure
theorem is a direct consequence of Theorems \ref{thm53}, \ref{thm64} and Propositions \ref{prop62}, \ref{prop63}.
\begin{theorem} \label{th71}
Let $A$ be a cycle-finite algebra. Then there exist quotient algebras $B_1, \ldots, B_p$ of $A$ such that the following statements hold.
\begin{enumerate}
\renewcommand{\labelenumi}{\rm(\roman{enumi})}
\item For each $i\in\{1,\ldots,p\}$, $B_i$ is either a tame generalized multicoil algebra or a tame generalized double tilted algebra.
\item $\ind A = \bigcup_{i=1}^p\ind B_i$.
\end{enumerate}
\end{theorem}

It follows from Theorem \ref{thm52} that for a tame tilted algebra $B$, all but finitely many modules in $\ind B$ are indecomposable modules over
the left tilted algebras $B_1^{(l)}, \ldots, B_m^{(l)}$ of Euclidean types or over the right tilted algebras $B_1^{(r)}, \ldots, B_n^{(r)}$
of Euclidean types. Moreover, all representation-infinite tilted algebras of Euclidean types and all tubular algebras are tame generalized
multicoil algebras. Therefore, applying Theorem \ref{thm43}, we obtain the following completion to Theorem \ref{th71}.
\begin{theorem} \label{th72}
Let $A$ be a cycle-finite algebra. Then there exist tame generalized multicoil quotient algebras $B_1, \ldots, B_q$ of $A$ whose indecomposable
finite dimensional modules exhaust all but finitely many isoclasses of modules in $\ind A$.
\end{theorem}

Moreover, we have the following information on nondirecting indecomposable modules over cycle-finite algebras.
\begin{theorem} \label{th73}
Let $A$ be a cycle-finite algebra. Then there exist tame generalized multicoil quotient algebras $B_1, \ldots, B_q$ of $A$ such that all but finitely many
isomorphism classes of nondirecting modules in $\ind A$ belong to generalized multicoils of the Auslander-Reiten quivers $\Gamma_{B_1}, \ldots, \Gamma_{B_q}$
of $B_1, \ldots, B_q$.
\end{theorem}

We exhibit an example of a cycle-finite algebra having a nongeneralized standard Auslander-Reiten component.
\begin{example}
Let $Q$ be the quiver
$$\xymatrix@C=16pt@R=12pt{
&2\ar[ld]_{\alpha}&3\ar[l]_{\gamma}\\
1&&&6\ar@<2pt>[lu]_{\eta}\ar@<-2pt>[lu]^{\xi}\\
&4\ar[lu]^{\beta}&5\ar[l]^{\sigma}\ar[ru]_{\delta}\\
}$$
$I$ the ideal in the path algebra $KQ$ of $Q$ generated by the elements $\gamma\alpha$, $\sigma\beta$, $\eta\gamma$, $\xi\gamma$, $\delta\xi$, $\delta\eta$,
and $A=KQ/I$ the associated bound quiver algebra. Then $A$ is a cycle-finite algebra whose Auslander-Reiten quiver $\Gamma_A$ has the disjoint union form
\[ \Gamma_A = {\mathcal C}\cup (\bigcup_{\lambda\in{\Bbb P}_1(K)}{\mathcal T}_{\lambda}), \]
where ${\mathcal T}_{\lambda}$, $\lambda\in{\Bbb P}_1(K)$, is the family of stable tubes of rank $1$ over the Kronecker path algebra $H=K\Delta$ given by
the subquiver $\Delta$ of $Q$ formed by the arrows $\xi$ and $\eta$, and $\mathcal C$ is the following glueing of the preprojective component ${\mathcal P}(H)$
and the preinjective component ${\mathcal Q}(H)$ of $\Gamma_A$ into a component of $\Gamma_A$
$$\xymatrix@C=10pt@R=10pt{
&&\tau_AI_3\ar@<2pt>[rd]\ar@<-2pt>[rd]&&I_3\ar@<2pt>[rd]\ar@<-2pt>[rd]\\
\cdots&{\phantom{I}}\ar@<2pt>[ru]\ar@<-2pt>[ru]&&\tau_AS_6\ar@<2pt>[ru]\ar@<-2pt>[ru]&&S_6\ar[rd]&&I_4\ar[rd]\\
&&&&&&P_5\ar[ru]\ar[rd]&&S_5\\
&&&P_2\ar[rd]&&S_4\ar[ru]&&I_6\ar[ru]\\
&&S_1\ar[ru]\ar[rd]&&I_1\ar[ru]\ar[rd]&&&&P_6\ar@<2pt>[rd]\ar@<-2pt>[rd]&&\tau_A^-P_6\ar@<2pt>[rd]\ar@<-2pt>[rd]\\
&&&P_4\ar[ru]&&S_2\ar[rd]&&S_3\ar@<2pt>[ru]\ar@<-2pt>[ru]&&\tau_A^-S_3\ar@<2pt>[ru]\ar@<-2pt>[ru]&&{\phantom{I}}&\cdots\\
&&&&&&P_3\ar[ru]\\
}$$
Observe that the simple module $S_3$ is the socle of the injective module $I_3$ and the canonical monomorphism $S_3\to I_3$ belongs to $\rad^{\infty}(\modu A)$,
because it factors through any stable tube ${\mathcal T}_{\lambda}$. Hence $\mathcal C$ is not generalized standard. On the other hand, every cycle in $\ind A$
consists of modules of a stable tube ${\mathcal T}_{\lambda}$ and hence is finite, because the stable tubes ${\mathcal T}_{\lambda}$ are generalized standard.
Therefore, $A$ is a cycle-finite algebra.
\end{example}

We end this section with an example of an algebra $A$ with partially ordered Auslander-Reiten components of $\Gamma_A$ which is not cycle-finite.
\begin{example}
Let $Q$ be the quiver
\[\xymatrix@C=16pt@R=16pt{
&6\ar[r]^{\xi}\ar[dd]_{\delta}&7\ar[d]^{\eta}\cr
1&&4\ar[ld]_{\gamma}\cr
&3\ar[lu]_{\alpha}\ar[ld]^{\beta}\cr
2&&5\ar[lu]^{\sigma}\cr
}\]
$I$ the ideal in the path algebra $KQ$ of $Q$ generated by the elements $\delta\alpha$, $\delta\beta$, $\xi\eta$ and $\eta\gamma$, and $A=KQ/I$ the
associated bound quiver algebra. Then $A$ is a tame algebra whose Auslander-Reiten quiver $\Gamma_A$ has the disjoint union form
\[ \Gamma_A = \mathcal{P}(A)\cup\mathcal{T}^A\cup\mathcal{C}, \]
where $\mathcal{P}(A)$ is a preprojective component, $\mathcal{T}^A=(\mathcal{T}_{\lambda}^A)_{\lambda\in\mathbb{P}_1(K)\setminus\{1\}}$
is a family of pairwise orthogonal generalized standard stable tubes, and $\mathcal{C}$ is a component of the form below
\begin{center}
\setlength{\unitlength}{0.7mm}
\begin{picture}(100,107)(-50,-100)

%% wierzcholki malego%%%%%%%%%%%%%%%

\put(0,0){\makebox(0,0){$\circ$}}
\put(-10,-10){\makebox(0,0){$\circ$}}
\put(-10,-10){\makebox(0,0){$\circ$}}
\put(-10,-20){\makebox(0,0){$\circ$}}
\put(-10,-30){\makebox(0,0){$\circ$}}
\put(-20,-20){\makebox(0,0){$\circ$}}
\put(-20,-40){\makebox(0,0){$\circ$}}
\put(-30,-10){\makebox(0,0){$\circ$}}
\put(-30,-20){\makebox(0,0){$\circ$}}
\put(-30,-30){\makebox(0,0){$\circ$}}
\put(-40,-20){\makebox(0,0){$\circ$}}
\put(-40,-40){\makebox(0,0){$\circ$}}
\put(-50,-10){\makebox(0,0){$\circ$}}
\put(-50,-20){\makebox(0,0){$\circ$}}
\put(-50,-30){\makebox(0,0){$\circ$}}
\put(-50,-40){\makebox(0,0){$\circ$}}
\put(-50,-50){\makebox(0,0){$\circ$}}
\put(-60,-20){\makebox(0,0){$\circ$}}
\put(-60,-40){\makebox(0,0){$\circ$}}
\put(-70,-10){\makebox(0,0){$\circ$}}
\put(-70,-20){\makebox(0,0){$\circ$}}
\put(-70,-30){\makebox(0,0){$\circ$}}
\put(-70,-40){\makebox(0,0){$\circ$}}
\put(-70,-50){\makebox(0,0){$\circ$}}

%%strzalki malego %%%%%%%%%%%%%%%%%%%%%%%%%%%

\multiput(-49,-49)(10,10){5}{\vector(1,1){8}}
\multiput(-69,-49)(10,10){4}{\vector(1,1){8}}
\multiput(-69,-29)(10,10){2}{\vector(1,1){8}}
\multiput(-19,-39)(10,10){1}{\vector(1,1){8}}

\multiput(-69,-31)(10,-10){2}{\vector(1,-1){8}}
\multiput(-69,-11)(10,-10){3}{\vector(1,-1){8}}
\multiput(-49,-11)(10,-10){3}{\vector(1,-1){8}}
\multiput(-29,-11)(10,-10){2}{\vector(1,-1){8}}

\multiput(-69,-20)(10,0){6}{\vector(1,0){8}}
\multiput(-69,-40)(10,0){3}{\vector(1,0){8}}

%% wiercholki duzego

\put(10,-10){\makebox(0,0){$\circ$}}
\put(30,-10){\makebox(0,0){$\circ$}}
\put(50,-10){\makebox(0,0){$\circ$}}
\put(70,-10){\makebox(0,0){$\circ$}}
\put(10,-30){\makebox(0,0){$\circ$}}
\put(30,-30){\makebox(0,0){$\circ$}}
\put(50,-30){\makebox(0,0){$\circ$}}
\put(70,-30){\makebox(0,0){$\circ$}}
\put(10,-50){\makebox(0,0){$\circ$}}
\put(30,-50){\makebox(0,0){$\circ$}}
\put(50,-50){\makebox(0,0){$\circ$}}
\put(70,-50){\makebox(0,0){$\circ$}}
\put(10,-70){\makebox(0,0){$\circ$}}
\put(30,-70){\makebox(0,0){$\circ$}}
\put(50,-70){\makebox(0,0){$\circ$}}
\put(70,-70){\makebox(0,0){$\circ$}}

\put(20,-20){\makebox(0,0){$\circ$}}
\put(20,-40){\makebox(0,0){$\circ$}}
\put(20,-60){\makebox(0,0){$\circ$}}
\put(40,-20){\makebox(0,0){$\circ$}}
\put(40,-40){\makebox(0,0){$\circ$}}
\put(40,-60){\makebox(0,0){$\circ$}}
\put(60,-20){\makebox(0,0){$\circ$}}
\put(60,-40){\makebox(0,0){$\circ$}}
\put(60,-60){\makebox(0,0){$\circ$}}

\put(10,-90){\makebox(0,0){$\circ$}}
\put(30,-90){\makebox(0,0){$\circ$}}
\put(50,-90){\makebox(0,0){$\circ$}}
\put(70,-90){\makebox(0,0){$\circ$}}
\put(20,-80){\makebox(0,0){$\circ$}}
\put(40,-80){\makebox(0,0){$\circ$}}
\put(60,-80){\makebox(0,0){$\circ$}}

%%%% strzalki duzego

\multiput(1,-1)(10,-10){7}{\vector(1,-1){8}}
\multiput(31,-11)(10,-10){4}{\vector(1,-1){8}}
\multiput(51,-11)(10,-10){2}{\vector(1,-1){8}}
\multiput(11,-31)(10,-10){6}{\vector(1,-1){8}}
\multiput(11,-51)(10,-10){4}{\vector(1,-1){8}}
\multiput(11,-71)(10,-10){2}{\vector(1,-1){8}}

\multiput(11,-29)(10,10){2}{\vector(1,1){8}}
\multiput(11,-49)(10,10){4}{\vector(1,1){8}}
\multiput(11,-69)(10,10){6}{\vector(1,1){8}}

\multiput(11,-89)(10,10){6}{\vector(1,1){8}}

\multiput(31,-89)(10,10){4}{\vector(1,1){8}}
\multiput(51,-89)(10,10){2}{\vector(1,1){8}}

%%linie przerywane

\multiput(10,-31)(0,-5){15}{\line(0,-1){2}}
\multiput(70,-11)(0,-5){19}{\line(0,-1){2}}

\multiput(-75,-20)(-3,0){3}{\makebox(0,0){$\cdot$}}
\multiput(-75,-40)(-3,0){3}{\makebox(0,0){$\cdot$}}

\multiput(30,-95)(0,-3){3}{.}
\multiput(50,-95)(0,-3){3}{.}

%% opis

\put(0,4){\makebox(0,0){$P_7$}}
\put(-6,-11){\makebox(0,0){$S_4$}}
\put(-6,-21){\makebox(0,0){$S_5$}}
\put(-6,-31){\makebox(0,0){$S_6$}}
\put(6,-11){\makebox(0,0){$S_7$}}
\put(-20,-15){\makebox(0,0){$I_3$}}
\put(-17,-42){\makebox(0,0){$I_7$}}
\put(-50,-44){\makebox(0,0){$I_1$}}
\put(-50,-54){\makebox(0,0){$I_2$}}
\put(14,-20){\makebox(0,0){$P_6$}}
\put(7,-31){\makebox(0,0){$S_3$}}
\put(70,-6){\makebox(0,0){$S_3$}}

\end{picture}
\end{center}
\vspace{0,7cm}
(see \cite[Example2.6]{[JMS]}). Then we obtain that the components of $\Gamma_A$ are partially ordered in $\ind A$. In particular, we conclude
that for every cycle $X_0\to X_1\to\cdots\to X_{r-1}\to X_r=X_0$ in $\ind A$, all modules $X_0, X_1, \ldots, X_{r-1}$ belong to the same component
of $\Gamma_A$. On the other hand, we have in $\ind A$ a cycle
\[ I_3 \buildrel {f_1}\over {\hbox to 8mm{\rightarrowfill}} P_7 \buildrel {f_2}\over {\hbox to 8mm{\rightarrowfill}} P_6 \buildrel {f_3}\over {\hbox to 8mm{\rightarrowfill}} I_3 \]
consisting of modules from $\mathcal C$ and with $f_3$ in $\rad^{\infty}(\modu A)$, because there is no path in $\mathcal C$ from $P_6$ to $I_3$.
Therefore, $A$ is not a cycle-finite algebra.
\end{example}

\section{Discrete indecomposable modules} \label{sec:8}

The aim of this section is to establish a common bound on the number of discrete indecomposable modules of a fixed dimension vector over
a cycle-finite algebra.
\begin{lemma} \label{lem81}
Let $A$ be an algebra, $\Gamma$ a generalized standard generalized multicoil of $\Gamma_A$, $n$ the rank of $K_0(A)$, and ${\bf d}$ a nonnegative
vector of $K_0(A)$.
Then the number of indecomposable modules $X$ in $\Gamma$ with $\vdim X$ $={\bf d}$ is bounded by $n$. Moreover, if $\Gamma$ consists of modules
which do not lie on infinite cycles then the number of indecomposable modules $X$ in $\Gamma$ with $\vdim X$ $={\bf d}$ is bounded by $n-1$.
\end{lemma}
\begin{proof}
Without loss of generality we may assume that $A$ is the support algebra of $\Gamma$. Let $\Gamma$ be an arbitrary generalized multicoil of $\Gamma_A$
which is generalized standard. We shall prove our claim by induction on the number $m$ of admissible operations which we have to do on a finite family
${\mathcal T}_1, {\mathcal T}_2, \ldots, {\mathcal T}_s$ of generalized standard stable tubes in order to obtain the generalized multicoil $\Gamma$.
If  $m=1$, then we can only do the admissible operation (ad~1) or (ad~1$^*$), so $s=1$. In this case our statement follows from
\cite[Proposition 4.1]{[Ma2]}.

Let $m>1$.
If the $m$th admissible operation is of type (ad~1), then by definition of (ad~1) we have $\vdim V = \vdim W$ for any modules $V, W\in\{Z_{ij},X'_i\}$,
$i\geq 0$, $1\leq j\leq t$ such that $V\neq W$.
%we know that the set of new vertices contains no two modules $V$, $W$ such that $\vdim V = \vdim W$.
Therefore, the number of indecomposable modules with the same dimension vector does not change.
If it is of type (ad~1$^*$) then the situation is the same.
If the $m$th admissible operation is of type (ad~2), then in the sequence of earlier $m-1$ admissible operations, there is an operation of type (ad~1$^*$)
or (ad~5$^*$) which contains an operation (fad~1$^*$) which gives rise to the pivot $X$ of (ad~2), and the operations done between these two must
not affect the support of $\Hom_A(X,-)$ restricted to the generalized multicoil containing $X$.
Let $t$ be the parameter of such operation (ad~1$^*$) or like in definition of (fad~1$^*$).
Note that in general, in the sequence of earlier $m-1$
admissible operations can be an operation of type (ad~5) which contains an operation (fad~4) which gives rise to the pivot $X$
of (ad~2) but from Lemma \cite[Lemma 3.10]{[MS2]} this case can be reduced to (ad~5$^*$) which contains an operation (fad~1$^*$).
Moreover, from \cite[Lemma 3.3]{[Ma2]} we know that, for a fixed ${\bf e}\in K_0(A)$, each ray and coray in $\Gamma$ contains at most one module $Z$
with $\vdim Z$ $={\bf e}$. Therefore, we get that the number of new indecomposable modules with the same dimension vector is at most $t+1$ while
after applying operations (ad~1$^*$) and (ad~2) or (fad~1$^*$) and (ad~2) the number of new vertices in ordinary quiver of $A$ increases by $t+2$.
If the $m$th admissible operation is of type (ad~2$^*$), then the proof is dual.
If the $m$th admissible operation is of type (ad~3), then in the sequence of earlier $m-1$ admissible operations, there must be at least one
operation of type (ad~1$^*$) or (ad~5$^*$) which contains the operation (fad~1$^*$) which gives rise to the pivot $X$ of (ad~3) and to the modules
$Y_1, Y_2, \ldots, Y_t$ in the support of $\Hom_A(X,-)$ restricted to the generalized multicoil containing $X$. The operations done after must not
affect this support. Again, in general, in the sequence of earlier $m-1$ admissible operations can be an operation of type (ad~5) which contains
an operation (fad~4) which gives rise to the pivot $X$ of (ad~3) but from \cite[Lemma 3.10]{[MS2]} this case can be reduced to (ad~5$^*$) which
contains an operation (fad~1$^*$).
Suppose that we had $r$ such consecutive operations of types (ad~1$^*$) or (fad~1$^*$), the first of which had $X_t$ as a pivot, and these admissible
operations built up a~branch $L$ in $A$ with points $a, a_1, a_2, \ldots, a_t$ in $Q_A$, so that $X_{t-1}$ and $Y_t$
are the indecomposable injective $A$-modules corresponding respectively to $a$ and $a_1$, and both $Y_1$ and $\tau_A^{-1}Y_1$ are coray modules
in the generalized multicoil containing the (ad~3)-pivot $X$.
Again, from \cite[Lemma 3.3]{[Ma2]} we know that, for a fixed ${\bf e}\in K_0(A)$, each ray and coray in $\Gamma$ contains at most one module $Z$
with $\vdim Z$ $={\bf e}$. Therefore, we get that the number of new indecomposable modules with the same dimension vector is at most $t+1$ while
after applying $r$ consecutive operations of types (ad~1$^*$) and (ad~3) or $r$ consecutive operations of types (fad~1$^*$) and (ad~3)
the number of new vertices in ordinary quiver of $A$ increases by $t+2$.
If the $m$th admissible operation is of type (ad~3$^*$), then the proof is dual.
If the $m$th admissible operation is of type (ad~4), then $A$ is the algebra obtained from another one, say $A'$, by applying this admissible
operation with pivot $X$ and the begin $Y_1$ of a finite sectional path $Y_1 \rightarrow Y_2 \rightarrow \cdots \rightarrow Y_t$.
Note that this finite sectional path is the linearly oriented quiver of type ${\Bbb A}_t$ and its support algebra $\Lambda$ (given by the vertices
corresponding to the simple composition factors of the modules $Y_1, Y_2, \ldots, Y_t$) is a tilted algebra of the path algebra $D$ of the linearly
oriented quiver of type ${\Bbb A}_t$. From \cite[(4.4)(2)]{[Ri]} we know that $\Lambda$ is a bound quiver algebra given by a truncated branch in $x$,
where $x$ corresponds to the unique projective-injective $D$-module. Moreover, the modules $Y_1, Y_2, \ldots, Y_t$ are directing in $\Gamma_{A'}$.
Again, from \cite[Lemma 3.3]{[Ma2]} we know that, for a fixed ${\bf e}\in K_0(A)$, each ray and coray in $\Gamma$ contains at most one module $Z$ with
$\vdim Z$ $={\bf e}$.
Therefore, we get that the number of new indecomposable modules with the same dimension vector is at most $t+r+1$ while
after applying operations which give rise to the finite sectional path $Y_1 \rightarrow Y_2 \rightarrow \cdots \rightarrow Y_t$ and (ad~4)
the number of new vertices in ordinary quiver of $A$ increases by $t+r+2$.
If the $m$th admissible operation is of type (ad~4$^*$), then the proof is dual.
There remains to consider the case where the $m$th admissible operation is of type (ad~5). Since in the definition of admissible operation (ad~5)
we use the finite versions (fad~1), (fad~2), (fad~3), (fad~4) of the admissible operations (ad~1), (ad~2), (ad~3), (ad~4) and the admissible operation
(ad~4), we conclude that the lemma follows from the above considerations. If it is of type (ad~5$^*$), then the proof is dual and this finishes the
proof of the lemma.
\qed
\end{proof}

The following proposition is essential for the proof of Theorem \ref{thm83}.

\begin{proposition} \label{prop82}
Let $A$ be a tame generalized multicoil algebra, $n$ the rank of $K_0(A)$, and ${\bf d}$ a nonnegative vector of $K_0(A)$.
Then
\begin{enumerate}
\renewcommand{\labelenumi}{\rm(\roman{enumi})}
\item The number of isomorphism classes of discrete indecomposable $A$-modules $X$ with $\vdim X$ $={\bf d}$ is bounded by $n+2$.
\item The number of isomorphism classes of indecomposable $A$-modules $X$ with $\vdim X$ $={\bf d}$ and $q_A(\vdim X)\neq 0$ is bounded by $n-1$.
\item The number of isomorphism classes of indecomposable $A$-modules $X$ with $\vdim X$ $={\bf d}$ and $\chi_A(\vdim X)\neq 0$ is bounded by $n-1$.
\end{enumerate}
\end{proposition}
\begin{proof}
It is known that an indecomposable $A$-module $M$ in $\Gamma_A$ which does not lie on oriented cycle in $\Gamma_A$ is uniquely determined by $[M]$.
Moreover, if $\mathcal T$ is a stable tube in $\Gamma_A$ then the support of $\mathcal T$ is a tame concealed or tubular convex subcategory of $A$.
Hence, for any indecomposable $A$-module $X$ lying in a stable tubes of rank one, we have $q_A(\vdim X)=\chi_A(\vdim X)=0$.
Let ${\bf d}$ be a nonnegative vector in $K_0(A)$ such that there exists a nondirecting, discrete indecomposable $A$-module $X$ with $\vdim X$ $={\bf d}$.
Then $X$ belongs to a generalized multicoil $\Gamma$ of $\Gamma_A$. Assume first that $\Hom_A(P,X)\neq 0$ for some indecomposable projective module
in $\Gamma$. Then it follows from the proof of \cite[Proposition 3.5]{[Ma2]} that any indecomposable $A$-module $Y$ with $\vdim Y$ $={\bf d}$
also lies in $\Gamma$. Applying now Lemma \ref{lem81} we conclude that the number of isomorphism classes of indecomposable $A$-modules $Z$ with
$\vdim Z = \vdim X$ $={\bf d}$ is bounded by $n-1$. We get the same statement in the case when $\Hom_A(X,I)\neq 0$ for an indecomposable injective
module $I$ in $\Gamma$.
Note that different tame concealed algebras and different tubular algebras give modules in $A$ with different dimension vectors.
Hence, it remains to consider the case when the support of $X$ is contained in a convex subcategory, say $C$, which is tame concealed or tubular.
Then $X$ belongs to a ${\Bbb P}_1(K)$-family ${\mathcal T}=({\mathcal T}_{\lambda})_{\lambda\in{\Bbb P}_1(K)}$ of standard stable tubes of
$\Gamma_C$. Moreover, if $Z$ is a indecomposable $A$-module with $\vdim Z = \vdim X$ $={\bf d}$ then $Z$ is a $C$-module and lies in one of the tubes
${\mathcal T}_{\lambda}$ (see \cite{[Ri]} or \cite{[Sk7]}). Denote by $m$ the rank of $K_0(C)$, and by $r_{\lambda}$ the rank of  the tube
${\mathcal T}_{\lambda}$, $\lambda\in{\Bbb P}_1(K)$. Then the following equality holds
\[ \sum_{\lambda\in {\Bbb P}_1(K)}(r_{\lambda}-1)=m-2 \]
(see \cite{[Sk7]}). Further, if $Y\in {\mathcal T}_{\lambda}$ and $Z\in {\mathcal T}_{\mu}$ are two nonisomorphic modules in $\mathcal T$ with
$\vdim Y = \vdim Z$ then the quasi-length of $Y$  is divisible by $r_{\lambda}$ and the quasi-length of $Z$ is divisible by $r_{\mu}$. We note that then
$q_A(\vdim Y)=q_C(\vdim Y)=\chi_A(\vdim Y)=0$ and $q_A(\vdim Z)=q_C(\vdim Z)=\chi_A(\vdim Z)=0$, since $\gldim C\leq 2$.
Now a simple inspection of tubular types of tame concealed and tubular algebras shows that, if $\lambda_1,\ldots,\lambda_t$ are all indices
$\lambda\in{\Bbb P}_1(K)$ with $r_{\lambda}\neq 1$, then $r_{\lambda_1}+\ldots+r_{\lambda_t}\leq m+2\leq n+2$.
Therefore, the number of isomorphism classes of indecomposable $A$-modules $Z$ with $\vdim Z = \vdim X$ $={\bf d}$ is bounded by $n+2$.
\qed
\end{proof}
\begin{theorem} \label{thm83}
Let $A$ be a cycle-finite algebra. Then there is a positive integer $m$ such that, for each nonnegative vector
${\bf d}\in K_0(A)$, the number of isomorphism classes of discrete indecomposable $A$-modules of dimension vector $\bf d$ is bounded by $m$.

Moreover, if $A$ is coherent, $n$ is the rank of $K_0(A)$, then the following statements hold:
\begin{enumerate}
\renewcommand{\labelenumi}{\rm(\roman{enumi})}
\item The number of isomorphism classes of discrete indecomposable $A$-modules $X$ with $\vdim X$ $={\bf d}$ is bounded by $n+2$.
\item The number of isomorphism classes of indecomposable $A$-modules $X$ with $\vdim X$ $={\bf d}$ and $q_A(\vdim X)\neq 0$ is bounded by $n-1$.
\item The number of isomorphism classes of indecomposable $A$-modules $X$ with $\vdim X$ $={\bf d}$ and $\chi_A(\vdim X)\neq 0$ is bounded by $n-1$.
\end{enumerate}
\end{theorem}
\begin{proof}
It follows from a result due to Ringel \cite[(2.4)(8)]{[Ri]} that, if $X, Y$  are modules in $\ind A$ with $\vdim X = \vdim Y$ and $X$
is directing, then $X$ and $Y$ are isomorphic. Further, by Theorem \ref{th73}, there exist tame generalized multicoil quotient algebras
$B_1, \ldots, B_q$ of $A$ such that all but finitely many isomorphism classes of nondirecting modules in $\ind A$ belong to generalized
multicoils of the Auslander-Reiten quivers $\Gamma_{B_1}, \ldots, \Gamma_{B_q}$ of $B_1, \ldots, B_q$. Moreover, it follows from the proof
of Theorem \ref{thm64} that, if $X$ and $Y$ are nondirecting discrete modules in $\ind A$ with $\vdim X = \vdim Y$ lying in generalized
multicoils of $\Gamma_{B_1}, \ldots, \Gamma_{B_q}$, then $X$ and $Y$ belong to the same generalized multicoil of $\Gamma_{B_p}$, for a fixed
$p\in\{1, \ldots, q\}$. Then there is a positive integer $m$ such that, for each nonnegative vector ${\bf d}\in K_0(A)$, the number of
isomorphism classes of discrete indecomposable $A$-modules of dimension ${\bf d}$ is bounded by $m$. Finally, if $A$ is coherent, then
all nondirecting modules in $\ind A$ lie in generalized multicoils of $\Gamma_{B_1}, \ldots, \Gamma_{B_q}$, and the statements (i), (ii), (iii)
are consequences of Proposition \ref{prop82}.
\qed
\end{proof}

\section{Homological properties of indecomposable modules} \label{sec:9}

Let $A$ be an algebra and $X$ be a nonprojective module in $\ind A$. Then we have an almost split sequence
\[ 0\to \tau_AX\to\bigoplus_{i=1}^{s(X)}Y_i\to X\to 0 \]
in $\modu A$ with $Y_1, \ldots, Y_{s(X)}$ indecomposable modules. We note that $\tau_A=D\Tr$, where $D=\Hom_K(-,K)$ is the standard duality on $\modu A$
and the transpose $\Tr X$ of a module $X$ in $\ind A$ is the cokernel of the homomorphism $\Hom_A(f,A)$ in $\modu A^{\op}$ associated to a minimal projective
presentation $P_1 \buildrel {f}\over {\hbox to 8mm{\rightarrowfill}} P_0\to X\to 0$ of $X$ in $\modu A$. Hence $\tau_A$ is a homological operator, and
$s(X)$ is a numerical homological invariant of a nonprojective module $X$ in $\ind A$.

We also recall that for an indecomposable nonsimple projective-injective module $P$ in $\modu A$ there is a canonical (up to isomorphism) almost split
sequence having $P$ as a middle term, namely
\[ 0\to \rad P \to P\oplus\rad P/\soc P \to P/\soc P \to 0 \]
(see \cite[Proposition IV.3.11]{[ASS]}).

The following theorem has been proved by Bautista and Brenner in \cite{[BaBr]} (see also \cite{[Li4]}).
\begin{theorem} \label{t91}
Let $A$ be an algebra of finite representation type and $X$ a nonprojective module in $\ind A$. Then
\begin{enumerate}
\renewcommand{\labelenumi}{\rm(\roman{enumi})}
\item $s(X)\leq 4$.
\item If $s(X)=4$, then one of the indecomposable middle terms $Y_i$ of an almost split sequence for $X$ is projective-injective.
\end{enumerate}
\end{theorem}

It has been conjectured by Brenner that, for any tame algebra $A$ and every nonprojective module $X$ in $\ind A$, we have $s(X)\leq 5$. It is still an open problem.
The following theorem proved in \cite[Theorem 3]{[PTa]} by the second named author and Takane confirms the Brenners's conjecture, and extends
Theorem \ref{t91}, to cycle-finite algebras.
\begin{theorem} \label{t92}
Let $A$ be a cycle-finite algebra and $X$ a nonprojective module in $\ind A$. Then
\begin{enumerate}
\renewcommand{\labelenumi}{\rm(\roman{enumi})}
\item $s(X)\leq 5$.
\item If $s(X)=5$, then one of the indecomposable middle terms $Y_i$ of an almost split sequence for $X$ is projective-injective.
\end{enumerate}
\end{theorem}

It has been proved by Ringel \cite[(2.4)(8)]{[Ri]} that, for any directing indecomposable module $X$ over an algebra $A$, we have $\End_A(X)\cong K$
and $\Ext_A^r(X,X)=0$ for $r\geq 0$. On the other hand, by \cite[Theorem B]{[MS3]}, for a tame generalized multicoil algebra $A$ and an arbitrary module $X$
in $\ind A$, we have $\dim_K\Ext_A^1(X,X)\leq \dim_K\End_A(X)$ and $\Ext_A^r(X,X)=0$ for $r\geq 2$. Hence, applying Theorems \ref{thm64} and \ref{thm65}
(and Corollary \ref{cor66}) we obtain the following theorem.
\begin{theorem}
Let $A$ be a cycle-finite algebra. Then for all but finitely many isomorphism classes of module $X$ in $\ind A$ we have $\dim_K\Ext_A^1(X,X)\leq \dim_K\End_A(X)$
and $\Ext_A^r(X,X)=0$ for $r\geq 2$.
\end{theorem}

We conclude from the above theorem that, for a cycle-finite algebra $A$ and all but finitely many isomorphism classes $X$ in $\ind A$, the Euler form
\[ \chi_A([X])=\sum_{r=0}^{\infty}(-1)^r\dim_K\Ext_A^r(X,X) \]
of $X$ is defined and is nonnegative. For $A$ a coherent cycle-finite algebra, it is the case for all modules $X$ in $\ind A$.

It is known from \cite[Theorem 2.3]{[HRS]} and \cite[Theorem 3.4]{[RS2]} that if $A$ is a quasitilted algebra or a generalized double tilted algebra then,
for all but finitely many isomorphism classes of modules $X$ in $\ind A$, we have $\pd_AX\leq 1$ or $\id_AX\leq 1$. Moreover, it has been conjectured in
\cite{[Sk11]} that the converse implication also holds.

We end this section by the following recent result by the third named author and Skowyrski \cite{[SkSk]} which confirms this conjecture for cycle-finite algebras.
\begin{theorem}
Let $A$ be a cycle-finite algebra such that $\pd_AX\leq 1$ or $\id_AX\leq 1$ for all but finitely many isomorphism classes of modules in $\ind A$.
Then $A$ is a tame quasitilted algebra or a tame generalized double tilted algebra.
\end{theorem}

\section{Geometric properties of indecomposable modules} \label{sec:10}

The aim of this section is to present some results describing geometric properties of indecomposable modules over cycle-finite algebras.

Let $A$ be an algebra and $A\cong KQ/I$ its bound presentation. Then $I$ is an admissible ideal in the path algebra $KQ$ of $Q$ generated by a finite
system of forms $\sum_{1\leq j\leq t}\lambda_j\alpha_{m_j,j}\ldots\alpha_{1,j}$ (called $K$-linear relations), where $\lambda_1, \ldots, \lambda_t$
are elements of $K$ and $\alpha_{m_j,j}, \ldots, \alpha_{1,j}$, $1\leq j\leq t$, are paths of length $\geq 2$ in $Q$ having a common source and a common
end. Denote by $Q_0$ the set of vertices of $Q$, by $Q_1$ the set of arrows of $Q$, and by $s, e: Q_1\to Q_0$ the maps which assign to each arrow
$\alpha_1$ its source $s(\alpha)$ and its end $e(\alpha)$. The category $\modu A$ of finite dimensional right $A$-modules is equivalent to the category
$\rep_K(Q,I)$ of all finite dimensional representations $V=(V_i,\varphi_{\alpha})_{i\in Q_0,\alpha\in Q_1}$ of $Q$, where $V_i$, $i\in Q_0$, are finite
dimensional $K$-vector spaces and $\varphi_{\alpha}: V_{s(\alpha)}\to V_{e(\alpha)}$, $\alpha\in Q_1$, are $K$-linear maps satisfying the equations
$\sum_{1\leq j\leq t}\lambda_j\varphi_{\alpha_{m_j,j}}\ldots\varphi_{\alpha_{1,j}}=0$ for all $K$-linear relations
$\sum_{1\leq j\leq t}\lambda_j\alpha_{m_j,j}\ldots\alpha_{1,j}\in I$ (see \cite{[ASS]}). Fix now a vector ${\bf d} = (d_i)_{i\in Q_0}\in K_0(A)={\Bbb Z}^{Q_0}$
with nonnegative coordinates. Denote by $\modu_A({\bf d})$ the set of all representations $V=(V_i,\varphi_{\alpha})$ in $\rep_K(Q,I)$ with $V_i=K^{d_i}$
for all $i\in Q_0$. A representation $V$ in $\modu_A({\bf d})$ is given by $d_{e(\alpha)}\times d_{s(\alpha)}$-matrices $V(\alpha)$ determining the maps
$\varphi_{\alpha}: K^{s(\alpha)}\to K^{e(\alpha)}$, $\alpha\in Q_1$, in the canonical bases of $K^{d_i}$, $i\in Q_0$. Moreover, the matrices $V(\alpha)$,
$\alpha\in Q_1$, satisfy the relations
\[ \sum_{1\leq j\leq t}\lambda_jV(\alpha_{m_j,j})\ldots V(\alpha_{1,j})=0 \]
for all $K$-linear relations $\sum_{1\leq j\leq t}\lambda_j\alpha_{m_j,j}\ldots\alpha_{1,j}\in I$. Therefore, $\modu_A({\bf d})$ is a closed subset of
${\Bbb A}({\bf d}) = \prod_{\alpha\in Q_1}K^{d_{e(\alpha)}\times d_{s(\alpha)}}$ in the Zariski topology, and so $\modu_A({\bf d})$ is an affine variety.
We note that $\modu_A({\bf d})$ is not necessarily irreducible. The affine (reductive) algebraic group $G({\bf d})=\prod_{i\in Q_0}\GL_{d_i}(K)$ acts
on the variety $\modu_A({\bf d})$ by conjugation
\[ (gV)(\alpha) = g_{e(\alpha)}V(\alpha)g_{s(\alpha)}^{-1} \]
for $g=(g_i)\in G({\bf d})$, $V\in\modu_A({\bf d})$, $\alpha\in Q_1$. We shall identify an $A$-module $V$ of dimension vector $\bf d$ with the
corresponding point of the variety $\modu_A({\bf d})$. The $G({\bf d})$-orbit $G({\bf d})M$ of a module $M$ in $\modu_A({\bf d})$ will be denoted by
$\mathcal{O}(M)$. Observe that two modules $M$ and $N$ in $\modu_A({\bf d})$ are isomorphic if an only if $\mathcal{O}(M)=\mathcal{O}(N)$.
For $M, N\in\modu_A({\bf d})$, we say that $N$ is a {\it degeneration} of $M$ if $N$ belongs to the Zariski closure $\overline{\mathcal{O}(M)}$ of $\mathcal{O}(M)$
in $\modu_A({\bf d})$, and we denote this fact by $M\leq_{\deg}N$. We note that $\leq_{\deg}$ is a partial order in $\modu_A({\bf d})$.
If $N\in\overline{\mathcal{O}(M)}$ implies $\mathcal{O}(N)=\mathcal{O}(M)$, then the orbit $\mathcal{O}(N)$ is said to be {\it maximal}.
Clearly, an orbit in $\modu_A({\bf d})$ of maximal dimension is maximal, but the converse is not true in general. It is known that the union of all
$G({\bf d})$-orbits in $\modu_A({\bf d})$ of maximal dimension is an open subset of $\modu_A({\bf d})$, called an {\it open sheet}
(see \cite{[Kr1]}, \cite{[Kr2]}). We note also that for a module $M$ in $\modu_A({\bf d})$, we have $\dim\mathcal{O}(M)=\dim G({\bf d})-\dim_K\End_A(M)$
(see \cite{[Kr1]}).
Given a module $M\in\modu_A({\bf d})$ we denote by $T_M(\modu_A({\bf d}))$ the tangent space of $\modu_A({\bf d})$ at $M$ and by $T_M(\mathcal{O}(M))$
the tangent space to $\mathcal{O}(M)$ at $M$. Then there is a canonical monomorphism of $K$-vector spaces
\[ T_M(\modu_A({\bf d}))/T_M(\mathcal{O}(M))\hookrightarrow \Ext_A^1(M,M) \]
(see \cite[(2.7)]{[Kr2]}. In particular, if $\Ext_A^1(M,M)=0$ then $\overline{\mathcal{O}(M)}$ is an irreducible component of $\modu_A({\bf d})$ and
$\mathcal{O}(M)$ is an open sheet of $\modu_A({\bf d})$. The {\it local dimension} $\dim_M\modu_A({\bf d})$ of $\modu_A({\bf d})$ is the maximal dimension
of the irreducible components of $\modu_A({\bf d})$ containing $M$. We have $\dim T_M(\modu_A({\bf d}))\geq \dim_M\modu_A({\bf d})$.
Further, $M\in\modu_A({\bf d})$ is said to be a {\it nonsingular point} of $\modu_A({\bf d})$ if $\dim_M\modu_A({\bf d})=\dim T_M(\modu_A({\bf d}))$.
If $M$ is a nonsingular point of $\modu_A({\bf d})$ then $M$ belongs to exactly one irreducible component of $\modu_A({\bf d})$ \cite[(II.2.6)]{[Sh]}.
The nonsingular points of $\modu_A({\bf d})$ form an open nonempty subset. It is known that a module $M$ in $\modu_A({\bf d})$ is nonsingular provided
$\Ext_A^2(M,M)=0$. A module variety $\modu_A({\bf d})$ is said to be a {\it complete intersection} provided the vanishing ideal of $\modu_A({\bf d})$ in the
coordinate ring $K[{\Bbb A}({\bf d})]$ of the affine space ${\Bbb A}({\bf d}) = \prod_{\alpha\in Q_1}K^{d_{e(\alpha)}\times d_{s(\alpha)}}$ is generated
by $\dim{\Bbb A}({\bf d})-\dim\modu_A({\bf d})$ polynomials. Finally, a module variety $\modu_A({\bf d})$ is said to be {\it normal} if the local ring
$\mathcal{O}_M$ of any module $M$ in $\modu_A({\bf d})$ is integrally closed in its total quotient ring. It is known that if $\modu_A({\bf d})$ is normal
then it is nonsingular in codimension one, that is, the set of singular points in $\modu_A({\bf d})$ is of codimension at most two
(see \cite[Chapter 11]{[E]}). If $\modu_A({\bf d})$ is a complete intersection, then $\modu_A({\bf d})$ is normal if and only if $\modu_A({\bf d})$
is nonsingular in codimension one (consequence of Serre's normality criterion). In the study of the degeneration order on a module variety
$\modu_A({\bf d})$ an important role is played by the following related partial orders. Let $M$ and $N$ be modules in a module variety. We define:
\begin{itemize}
\item $M\leq _{\ext} N$: $\Leftrightarrow$ there are modules $M_i$, $U_i$, $V_i$ and short exact sequences $0\to U_i\to M_i\to V_i\to 0$ in $\modu A$
such that $M=M_1$, $M_{i+1}=U_i\oplus V_i$, $1\leq i\leq s$, and $N=M_{s+1}$ for some natural number $s$.
\item $M\leq_R N$: $\Leftrightarrow$ there exists in $\modu A$ an exact sequence of the form $0\to N\to M\oplus Z\to Z\to 0$.
\item $M\leq N$: $\Leftrightarrow$ $\dim_K\Hom_A(M,X)\leq \dim_K\Hom_A(N,X)$ for all modules $X$ in $\modu A$.
\end{itemize}
It follows from the result due to Auslander \cite{[Au]} that $\leq$ is a partial order on the isomorphism classes of modules with the same dimension
vector. Further, for modules $M$ and $N$ in $\modu_A({\bf d})$, we have $M\leq N$ if and only if $\dim_K\Hom_A(X,M)\leq \dim_K\Hom_A(X,N)$
for all modules $X$ in $\modu A$ by a result of Auslander and Reiten \cite{[AR]}. Moreover, by a result of Zwara \cite{[Z1]}, we have
$M\leq_R N$ if and only if there exists in $\modu A$ a short exact sequence of the form $0\to Z'\to Z'\oplus M\to N\to 0$.

The following fundamental result of Zwara from \cite{[Z3]} (see also \cite{[Rie]} for the sufficiency part)
gives an algebraic characterization of degenerations of modules.
\begin{theorem}
Let $A$ be an algebra, $\bf d$ a vector in $K_0(A)$ with nonnegative coordinates, and $M$, $N$ modules in $\modu_A({\bf d})$. Then $M\leq_{\deg}N$
if and only if $M\leq_R N$.
\end{theorem}
In general, we have the following relations between the introduced orders. For modules $M$ and $N$ in the module variety $\modu_A({\bf d})$ the
following implications hold
\[ M\leq_{\ext} N \Longrightarrow M\leq_{\deg} N\Longrightarrow M\leq N. \]
Unfortunately, the reverse implications are not true in general, and it would be interesting to find out when there are true.

The following result of Zwara from \cite{[Z2]} gives a combinatorial description of degenerations for modules over algebras of finite
representation type.
\begin{theorem}
Let $A$ be an algebra of finite representation type, $\bf d$ a vector in $K_0(A)$ with nonnegative coordinates, and $M$, $N$ modules in $\modu_A({\bf d})$.
Then $M\leq_{\deg}N$ if and only if $M\leq N$.
\end{theorem}
We also exhibit the following results from \cite{[SZ2]} and \cite{[SZ3]} on degenerations of modules from the additive categories of generalized standard
Auslander-Reiten components.
\begin{theorem}
Let $A$ be an algebra, $\mathcal C$ a generalized standard quasi-tube of $\Gamma_A$, and $M$, $N$ modules in $\add({\mathcal C})$. Then $M\leq_{\deg}N$
if and only if $M\leq_{\ext} N$.
\end{theorem}
\begin{theorem}
Let $A$ be an algebra, $\mathcal C$ a generalized standard component of $\Gamma_A$, $N$ a module in $\add({\mathcal C})$, and $M$ a module in $\modu A$.
If $M\leq_{\deg}N$ then $M$ belongs to $\add({\mathcal C})$.
\end{theorem}
\begin{theorem}
Let $A$ be an algebra, $\mathcal C$ a generalized standard component of $\Gamma_A$, $M, N$ modules in $\add({\mathcal C})$ with $\vdim M = \vdim N$.
The following conditions are equivalent.
\begin{enumerate}
\renewcommand{\labelenumi}{\rm(\roman{enumi})}
\item $M\leq_{\deg}N$.
\item There exists an exact sequence $0\to N\to M\oplus Z\to Z\to 0$ in $\modu A$ with $Z$ from $\add({\mathcal C})$.
\item There exists an exact sequence $0\to Z'\to Z'\oplus M\to N\to 0$ in $\modu A$ with $Z'$ from $\add({\mathcal C})$.
\item $\dim_K\Hom_A(M,X)\leq \dim_K\Hom_A(N,X)$ for all modules $X$ in $\mathcal C$.
\item $\dim_K\Hom_A(X,M)\leq \dim_K\Hom_A(X,N)$ for all modules $X$ in $\mathcal C$.
\end{enumerate}
\end{theorem}

Let $A$ be an algebra and $M, N$ be nonisomorphic modules in $\ind A$ with $\vdim M = \vdim N$. Then $M\leq N$ forces the inequalities
$\dim_K\Hom_A(M,M)\leq \dim_K\Hom_A(N,M)$ and $\dim_K\Hom_A(M,M)\leq \dim_K\Hom_A(M,N)$, and consequently we have a cycle $M\to N\to M$.
Since $M\leq_{\deg}N$ implies $M\leq N$, we conclude that the directing modules in $\ind A$ are never involved in proper degenerations of
indecomposable modules. Observe also that, if $A$ is a cycle-finite algebra and $M<_{\deg}N$, then $M$ and $N$ belong to the same cyclic
component of $\Gamma_A$. The degenerations of modules in the additive categories of generalized multicoils of Auslander-Reiten quivers
of algebras were investigated in \cite{[Ma3]}, \cite{[Ma4]}, \cite{[SZ1]}, \cite{[SZ2]}. Then, using Theorems \ref{thm64} and \ref{thm65},
we obtain the following results.
\begin{theorem}
Let $A$ be a cycle-finite algebra. Then there exists a positive integer $t$ such that for any sequence
\[ M_r <_{\deg} M_{r-1} <_{\deg} \ldots <_{\deg} M_2 <_{\deg} M_1 \]
with $M_1, \ldots, M_r$ modules in $\ind A$, the inequality $r\leq t$ holds.
\end{theorem}
\begin{theorem}
Let $A$ be a coherent cycle-finite algebra and $M, M', N$ be modules in $\modu A$ such that $M <_{\deg}N$, $M' <_{\deg}N$ and $N$ is
indecomposable. Then $M\cong M'$ and is indecomposable.
\end{theorem}

The geometry of directing modules over tame algebras has been described in \cite{[BoS1]}. In particular, we have the following consequence
of \cite[Theorems 1 and 2]{[BoS1]}.
\begin{theorem}
Let $A$ be a tame algebra, $M$ a directing module in $\ind A$, and ${\bf d} = \vdim M$. Then the following statements hold.
\begin{enumerate}
\renewcommand{\labelenumi}{\rm(\roman{enumi})}
\item $\modu_A({\bf d})$ is a complete intersection and has at most two irreducible components.
\item The maximal $G({\bf d})$-orbits in $\modu_A({\bf d})$ consist of nonsingular modules.
\item $\mathcal{O}(M)$ is an open sheet of $\modu_A({\bf d})$.
\item All but finite number of $G({\bf d})$-orbits in $\modu_A({\bf d})$ have codimension at least two.
\item All $G({\bf d})$-orbits of codimension one are contained in $\overline{\mathcal{O}(M)}$.
\item If $M$ is not projective-injective over $\supp M$ then $\modu_A({\bf d})=\mathcal{O}(M)$, is normal and a complete intersection.
\end{enumerate}
\end{theorem}

It follows from Theorems \ref{thm31} and \ref{thm33} that the regular components of the Auslander-Reiten quivers of cycle-finite algebras
are stable tubes and their supports are tame concealed or tubular algebras. Then the following result from \cite[Theorem 1]{[BoS3]} describes
the geometry of modules from the additive categories of regular components of cycle-finite algebras.
\begin{theorem}
Let $A$ be a tame concealed or tubular algebra and $\bf d$ the dimension vector of a module in $\modu A$ which is periodic with respect
to the action of $\tau_A$. Then the affine variety $\modu_A({\bf d})$ is irreducible, normal, a complete intersection,
$\dim\modu_A({\bf d})=\dim G({\bf d}) - q_A({\bf d})$ and the maximal $G({\bf d})$-orbits in $\modu_A({\bf d})$ form the open sheet consisting
of nonsingular points. Moreover, a module $M$ in $\modu_A({\bf d})$ is a nonsingular point if and only if $\Ext_A^2(M,M) = 0$.
\end{theorem}

We refer also to \cite{[BoS2]} for the geometry of indecomposable modules over tame quasi-tilted algebras, which give information on the geometry
of indecomposable nondirecting modules of semiregular tubes of cycle-finite algebras (see Theorems \ref{thm36} and \ref{thm37}).

We end this section with the following consequence of \cite[Theorems A and B]{[MS3]} and Theorem \ref{thm64}, extending \cite[Theorem A]{[PS1]}
from strongly simply connected algebras of polynomial growth to coherent cycle-finite algebras.
\begin{theorem}
Let $A$ be a coherent cycle-finite algebra, $M$ a module in $\ind A$ and ${\bf d}=\vdim M$. Then the following statements hold.
\begin{enumerate}
\renewcommand{\labelenumi}{\rm(\roman{enumi})}
\item $M$ is a nonsingular point of $\modu_A({\bf d})$.
\item $q_A({\bf d})\geq\chi_A({\bf d})=\dim_K\End_A(M)-\dim_K\Ext_A^1(M,M)\geq 0$.
\item $\dim_M\modu_A({\bf d})=\dim G({\bf d})-\chi_A({\bf d})$.
\end{enumerate}
\end{theorem}

We note that there are indecomposable modules $M$ over coherent cycle-finite algebras $A$ with arbitrary large $\chi_A(\vdim M)$
(see \cite[(5.3)]{[PS2]}).

\begin{acknowledgement}
The authors gratefully acknowledge supports from the grant Maestro of the National Science Center of Poland and
the Centro de Investigaci\'on en Mathem\'aticas (CIMAT) Guanajuato in M\'exico.
\end{acknowledgement}

\end{document}